\documentclass[%
 aip,
 amsmath,amssymb,
preprint,%
]{revtex4-1}

\usepackage{graphicx}% Include figure files
\usepackage{dcolumn}% Align table columns on decimal point
\usepackage{bm}% bold math
\usepackage[utf8]{inputenc}
\usepackage[T1]{fontenc}
\usepackage{mathptmx}
\usepackage{etoolbox}
\usepackage{multirow}
\usepackage{amssymb}
\usepackage{color}
\usepackage{verbatim}
\usepackage{hyperref}
\usepackage{subfigure}
\usepackage{caption}
\usepackage{booktabs}
\usepackage{gensymb}
\usepackage{float}
\usepackage{amsmath}
\usepackage[version=4]{mhchem}
\usepackage{siunitx}
\usepackage{tabularx}
\makeatletter
\def\@email#1#2{%
 \endgroup
 \patchcmd{\titleblock@produce}
  {\frontmatter@RRAPformat}
  {\frontmatter@RRAPformat{\produce@RRAP{*#1\href{mailto:#2}{#2}}}\frontmatter@RRAPformat}
  {}{}
}%
\makeatother
\begin{document}

\preprint{AIP/123-QED}

\title{A component-splitting implicit time integration for multicomponent reacting flows simulations}
% Force line breaks with \\
\author{Jingchao Zhang}
\affiliation{School of Aeronautics, Northwestern Polytechnical University, Xi'an 710072, China}%Lines break automatically or can be forced with \\
\author{Jinsheng Cai}%
\affiliation{School of Aeronautics, Northwestern Polytechnical University, Xi'an 710072, China}

\author{Shucheng Pan}
\email{shucheng.pan@nwpu.edu.cn}
%  \homepage{http://www.Second.institution.edu/~Charlie.Author.}
\affiliation{School of Aeronautics, Northwestern Polytechnical University, Xi'an 710072, China}
% \affiliation{%
% Second institution and/or address%\\This line break forced% with \\
% }%

\date{\today}% It is always \today, today,
             %  but any date may be explicitly specified

	% \author[1]{Jingchao Zhang}
	% %\ead{zhangjingchao@mail.nwpu.edu.cn}
	% \ead{shucheng.pan@nwpu.edu.cn}
	% \cortext[cor1]{Corresponding author}
	% \author[1]{Jinsheng Cai}
	% \author[1,2]{Shucheng Pan\corref{cor1}}
	% \address[1]{School of Aeronautics, Northwestern Polytechnical University, Xi'an 710072, China}
	% \address[2]{Institute of Extreme Mechanics, Northwestern Polytechnical University, Xi'an 710072, China}

	%\affiliation{organization={NWPU},%Department and Organization
	%            addressline={}, 
	%            city={},
	%            postcode={}, 
	%            state={},
	%            country={}}

	\begin{abstract}
		A component-splitting method is proposed to improve convergence characteristics for implicit time integration of compressible multicomponent reactive flows. The characteristic decomposition of flux jacobian of multicomponent Navier-Stokes equations yields a large sparse eigensystem, presenting challenges of slow convergence and high computational costs for implicit methods. To addresses this issue, the component-splitting method segregates the implicit operator into two parts: one for the flow equations (density/momentum/energy) and the other for the component equations. Each part's implicit operator employs flux-vector splitting based on their respective spectral radii to achieve accelerated convergence. This approach improves the computational efficiency of implicit iteration, mitigating the quadratic increase in time cost with the number of species. Two consistence corrections are developed to reduce the introduced component-splitting error and ensure the numerical consistency of mass fraction. Importantly, the impact of component-splitting method on accuracy is minimal as the residual approaches convergence. The accuracy, efficiency, and robustness of component-splitting method are thoroughly investigated and compared with the coupled implicit scheme through several numerical cases involving thermo-chemical nonequilibrium hypersonic flows. The results demonstrate that the component-splitting method decreases the required number of iteration steps for convergence of residual and wall heat flux, decreases the computation time per iteration step, and diminishes the residual to lower magnitude. The acceleration efficiency is enhanced with increases in CFL number and number of species.% As a result, the component-splitting method improves convergence characteristics for implicit time integration of multicomponent Navier-Stokes equations.The first consistence correction ensures conservation properties by normalizing the iterative increment of the conservative variables. The second consistence correction normalizes the mass fraction which compromises the conservation property but enhances robustness. Despite the compromise in conservation,due to the coupled flux-vector splitting to approximate flux jacobian of implicit operator
	\end{abstract}

	%%Graphical abstract
	%\begin{graphicalabstract}
		%\includegraphics{grabs}
	%\end{graphicalabstract}

	%Research highlights
	% \begin{highlights}
	% 	\item A component-splitting implicit method is proposed for convergence acceleration for compressible multicomponent reacting flows.
	% 	\item Two type consistency corrections are developed to reduce splitting error and ensure numerical consistencies of mass fraction.
	% 	\item This method achieve better convergence characteristics on residual and wall heat flux for thermochemical nonequilibrium flows.
	% 	\item Large number of species and large CFL number shows better convergence acceleration effects of component-splitting method.  
	% \end{highlights}

	% \begin{keyword}
	% 	Convergence acceleration \sep multicomponent Navier-Stokes equations \sep implicit time integration \sep component-splitting \sep hypersonic flows
	% 	%% keywords here, in the form: keyword \sep keyword

	% 	%% PACS codes here, in the form: \PACS code \sep code

	% 	%% MSC codes here, in the form: \MSC code \sep code
	% 	%% or \MSC[2008] code \sep code (2000 is the default)

	% \end{keyword}

%% \linenumbers
\maketitle

% \begin{quotation}
% The ``lead paragraph'' is encapsulated with the \LaTeX\ 
% \verb+quotation+ environment and is formatted as a single paragraph before the first section heading. 
% (The \verb+quotation+ environment reverts to its usual meaning after the first sectioning command.) 
% Note that numbered references are allowed in the lead paragraph.
% %
% The lead paragraph will only be found in an article being prepared for the journal \textit{Chaos}.
% \end{quotation}
% \section{\label{sec:level1}First-level heading:\protect\\ The line
% break was forced \lowercase{via} \textbackslash\textbackslash}
\section{Introduction}
% 引言1、化学反应流动的应用2、隐式迭代的优势3、迭代速度较慢-波

In hypersonic flows, the intense shock compression and viscous dissipation lead to a substantial temperature rise, triggering chemical reactions within gas mixtures and vibrational energy excitation within molecules. Due to the comparable timescales of chemical reactions and hypersonic flows, fluid particles fail to reach a thermochemical equilibrium state before transitioning to the next position without undergoing sufficient collisions~\cite{anderson2000hypersonic}. Consequently, representing internal energy solely through a unified temperature becomes inadequate, necessitating the adoption of a multiple temperature model~\cite{park1989assessment}. As a result, hypersonic flows typically remains in a state of thermochemical non-equilibrium, and their numerical methods solving compressible multicomponent Navier-Stokes equations. %This multicomponent reactive system is also extensively applied in numerical simulations of combustion flows, describing multicomponent mixtures and detailed chemical reaction mechanisms over broad range of time scales due to the timescale of combustion may differ orders of magnitude from the timescale of flow.

The multicomponent Navier-Stokes equations entail conservation equations of vibrational-electronic energy and species mass, leading to a notable increase in computational consumption compared to thermo-chemical equilibrium gases. Additionally, the disparity in time scales within flows introduces considerable numerical stiffness to the solution, exacerbating challenges related to time integration. Time integration, which involves advancing the current solution to the next time step, is typically classified into explicit, implicit, and implicit-explicit methods. Explicit integration, exemplified by backward Euler~\cite{SHEN2014432} and multi-stage Runge-Kutta scheme~\cite{jameson1985numerical}, offer high-order accuracy but are constrained by stability requirements dictated by the Courant-Friedrichs-Lewy (CFL) condition. Implicit integration addresses this limitations by solving the conservation equations at the implicit time level, linearizing the equations through Taylor expansion on the right-hand side at the known time level. This method provides larger time steps and achieve superlinear convergence capabilities. In the context of unsteady flow, the implicit integration can be effectively combined with the dual time-stepping~\cite{jameson1991time} and local time-stepping~\cite{blazek2015computational} to attain favorable convergence rates while maintaining temporal discretization accuracy. Implicit-explicit time integration~\cite{doi:10.1137/0732037,LI2013157} have been developed to leverage the strengths of both explicit and implicit methods, such as employing implicit Runge-Kutta schemes~\cite{ASCHER1997151} to address stiff chemistry. 

The upwind spatial discretization schemes are commonly used to decompose convective fluxes into positive and negative fluxes based on the eigenvalues of the flux Jacobian matrix. Flux-difference splitting (FDS)~\cite{ROE1997250} and flux-vector splitting (FVS)~\cite{STEGER1981263} yields different approximations of the flux jacobian, including Roe splitting, eigenvalue splitting, and spectral radius (maximum eigenvalue) splitting. Implicit integration necessitates solving block-sparse systems, for which two main methods exist: direct methods, such as direct inversion via Gauss elimination, and iterative methods, including relaxation iterative algorithms~\cite{doi:10.2514/6.2006-2824}, approximation factorization decomposition\cite{BRILEY1977372}, and Newton-Krylov subspace methods\cite{doi:10.1137/0907058}. While direct methods can be computationally expensive for large systems, iterative methods are commonly employed for large-scale block-sparse systems. Point relaxation and line relaxation performs symmetric Gauss-Seidel iterations in both forward and backward directions using weighted average of neighboring values, with the weighting factor determined by the specific relaxation scheme employed~\cite{doi:10.1146/annurev.fl.07.010175.000431}. However, relaxation algorithms converge slowly or even diverge for ill-conditioned and non-diagonally dominant matrices on left-hand side. Factorization decomposition approximates the left-hand side matrix by a product of serval easily invertible matrix factors, such as lower-upper symmetric Gauss-Seidel (LUSGS)~\cite{doi:10.2514/3.10007} and alternating direction implicit (ADI)~\cite{BEAM197687}. Approximate factorization methods are more robust but introduce decomposition errors that can slow the convergence rate. Newton-Krylov subspace methods are efficient in solving large scale sparse linear systems~\cite{KNOLL2004357}, with one of the most widely used applications being the generalized minimum residual (GMRES)\cite{doi:10.1137/0907058} method. Newton-Krylov Subspace methods necessitate an initial solution close to the converged solution, thus requiring preconditioning. The convergence acceleration of implicit integration has garnered continuous attention from researchers. Jin Yao et al.~\cite{JIN2021110009} proposed a implicit boundary condition to introduce the flux jacobian from neighbouring multiblock for accelerated convergence of subiterative DDADI/D3ADI iteration. Zhang Rui et al.\cite{10.1063/5.0186520} developed conservative implicit scheme for solving the three-dimensional steady flows of molecular gases in all flow regimes from continuum one to free-molecular one. Bohao Zhou et al.~\cite{10.1063/5.0107571}proposed data-parallel upper-lower relaxation scheme based on Jacobi iteration to avoid dimensionality reduction in parallel computing, resulting in a higher convergence speed. Cao Wenbo et al.\cite{10.1063/5.0138863} proposed an online dimension reduction optimization method to enhance the convergence of the traditional iterative method.% Cavalca et al.\cite{CAVALCA2018399} developed a implicit defect-correction method for solving block-sparse system, demonstrating excellent efficiency and robustness for the Euler equations.

The implicit time integration encounters substantial challenges when solving multicomponent Navier-Stokes equations. One challenge from the disparity in timescale between chemical reactions and flows, while another arises from the difficulty in solving large-scale sparse matrices, impacting the convergence characteristics of implicit time integration. To enhance the computational efficiency, Gilbert Strang \cite{doi:10.1137/0705041} proposed a decoupled method that separates the chemical source terms from flow equations by solving stiff ODEs, known as the Strang splitting. Dong Haibo et al.~\cite{DONG2018146} developed a improved decoupled method which removes the standard enthalpy of formation from the total energy to thoroughly separate the effect of chemical reaction from flow equations, achieving higher computational efficiency for reactive flows. Hansen and Sutherland~\cite{doi:10.1137/15M1023166} developed a robust, efficient, and adaptive dual time time-stepping technique to alleviate the time step constraint associated with ignition/extinction in the context of combustion chemistry. Xianliang Chen and Song Fu~\cite{CHEN2020104668} investigated the preconditioner of GMRES in accelerating implicit time-stepping for hypersonic thermochemical nonequilibrium flows. Although the decoupled method can enhance the time step, the implicit time integration in solving multicomponent equations still face challenges in computational time and convergence. The computational complexity becomes prohibitive as the matrix dimension increases quadratically with the number of species. Wang Jian-Hang et al. ~\cite{WANG2019364} proposed a partial characteristic decomposition to improve the computational efficiency without deteriorating high-order accuracy for high-order finite difference schemes of multi-species euler equations by splitting the multicomponent eigen system into gas mixture part and species mass part. While this method addresses the issue of an exponential increase in computational cost for explicit time integration, the problem of significant computational consumption for a large number of species persists in implicit methods. %Ceze and Fidkowski~\cite{https://doi.org/10.1002/nme.4858} proposed a constrained pseudo transient continuation with positivity preservation for an adaptive, high-order discontinuous Galerkin method for Reynolds-averaged Navier-Stokes simulation.

Aiming at accelerating the convergence of multicomponent Navier-Stokes equations using implicit time iteration, a component-splitting method is proposed to perform implicit iterations separately for the component conservation equations and the momentum/energy conservation equations. The left-hand side implicit operator are reformulated according to these two types of equations. The accuracy, efficiency, and robustness of this method are investigated through numerical simulations of hypersonic flows in a state of thermochemical non-equilibrium. %Furthermore, the component-splitting method also enhances the magnitude of residual reduction, which implies obtaining a more accurate solution by the implicit time integration.

\section{Methodology}
% \subsection{\label{sec:level2}Second-level heading: Formatting}
\subsection{Governing equations}
\begin{comment}

\begin{equation}
	\label{Eq.flux}
  \begin{array}{l}
  \bm{F} = u_1 \bm{Q} + P[0,1,0,0,u,0,\dots,0]^{\text{T}},\\
  \bm{G} = u_2 \bm{Q} + P[0,0,1,0,u_2,0,\dots,0]^{\text{T}}, \\
  \bm{H} = u_3 \bm{Q} + P[0,0,0,1,u_3,0,\dots,0]^{\text{T}}, \\
\end{array}
\end{equation}
and
\begin{equation}
	\label{Eq.flux}
  \begin{array}{l}
  \bm{F_v} = [0, \tau_{11}, \tau_{12}, \tau_{13}, u_i \tau_{i1} - q_{1}, - \rho Y_1 V_{1,1}, \dots, -\rho Y_{ns}V_{ns,1}]^{\text{T}},\\
  \bm{G_v} = [0, \tau_{21}, \tau_{22}, \tau_{23}, u_i \tau_{i2} - q_{2}, - \rho Y_1 V_{1,2}, \dots, -\rho Y_{ns}V_{ns,2}]^{\text{T}},\\
  \bm{H_v} = [0, \tau_{31}, \tau_{32}, \tau_{33}, u_i \tau_{i3} - q_{3}, - \rho Y_1 V_{1,3}, \dots, -\rho Y_{ns}V_{ns,3}]^{\text{T}},\\
\end{array}
\end{equation}
\end{comment}

The governing equations for the $l$-dimensional compressible multicomponent reactive Navier-Stokes systems are written as,

\begin{equation}
	\label{Eq:govern}
	\frac{{\partial {\bm{Q}}}}{{\partial t}} + \frac{{\partial {(\bm{F}_{inv}+\bm{F}_{vis})}}}{{\partial x_l}} = {\bm{S}},
\end{equation}
where $\bm{Q},\bm{F}_{inv}, \bm{F}_{vis}$, and $\bm{S}$ are the vector of conservative variables, inviscid (convective) fluxes, viscid fluxes, and source term respectively, i.e.,
\begin{equation}
	\label{Eq.flux}
	\begin{aligned}
	{\bm{Q}} = &\left[\rho,\rho u_1  ,\rho u_2  ,\rho u_3  ,\rho e_t ,\rho e_{v} ,\rho Y_1  ,\cdot ,\rho Y_{ns-1} \right]^{\mathrm{T}}\\
	 {\bm{F}_{inv}} = &  
	\left[\rho u_l,\rho u_1 u_l + p \delta_{1,l},\rho u_2 u_l + p \delta_{2,l},\rho u_3 u_l + p \delta_{3,l},\rho h u_l ,\rho e_v u_l,\right. &\\ 
	&\left. \rho Y_1 u_l,\cdot,\rho Y_{ns-1} u_l \right]^{\mathrm{T}} & \\
	{{\bm{F}}_{vis}} = & \left[0,- \tau_{1,l},-\tau_{2,l},-\tau_{3,l},-u_i \tau_{il} + q_{l},q_{v,l},\rho Y_1 V_{1,l},\cdot,\rho Y_{ns-1} V_{ns-1,l}\right]^{\mathrm{T}}\\
	{{\bm{S}}} = &\left[0,0,0 ,0,0,\omega_{v},\omega_{1},\cdot,\omega_{ns-1}\right]^{\mathrm{T}},
	\end{aligned}
\end{equation}
% \begin{equation}
% 	\label{Eq.flux}
% 	{\bm{Q}} = \left[ \begin{array}{l}
%     \rho                    \\
% 		\rho u_1  \\
% 		\rho u_2  \\
%     \rho u_3  \\
% 		\rho e_t \\
% 		\rho e_{v} \\
%     \rho Y_1  \\
%     \cdot \\
%     \rho Y_{ns-1} \\
% 	\end{array} \right], {\bm{F}_{inv}} =  
% 	\left[ \begin{array}{l}
% 		\rho u_l                    \\
%     \rho u_1 u_l + p \delta_{1,l} \\
%     \rho u_2 u_l + p \delta_{2,l} \\
%     \rho u_3 u_l + p \delta_{3,l}\\
% 		\rho h u_l \\
%     \rho e_v u_l \\
% 		\rho Y_1 u_l                              \\
%     \cdot \\
%     \rho Y_{ns-1} u_l                              \\
% 	\end{array} \right]
% 	,
% 	{{\bm{F}}_{vis}} = \left[ \begin{array}{l}
%     0 \\
%     - \tau_{1,l} \\
%     -\tau_{2,l} \\ 
%     -\tau_{3,l} \\
%     -u_i \tau_{il} + q_{l} \\
%     q_{v,l}\\
% 		\rho Y_1 V_{1,l}\\
%     \cdot \\
% 		\rho Y_{ns-1} V_{ns-1,l}\\
% 	\end{array} \right],
% 	{{\bm{S}}} = \left[ \begin{array}{l}
%     0 \\
%     0 \\
%     0 \\ 
%     0 \\
%     0 \\
%     \omega_{v}\\
% 		\omega_{1}\\
%     \cdot \\
% 		\omega_{ns-1}\\
% 	\end{array} \right],
% \end{equation}
where ${\rho}, e_t, e_v, h, \omega_v$ are the density, specific total energy, specific vibrational-electronic energy , specific enthalpy, vibrational energy source, respectively. $Y_{s}$ and $\omega_{s}$ are the mass fraction and chemical production rate of $s$th species, respectively. ${p}$ is the pressure, $\delta$ is the kronecker delta. $ns$ is the total number of chemical species. The density of $ns$th species is computed according to the conservation of mass, $Y_{ns} = 1 - \sum_{s=1}^{ns-1}Y_{s}$. $\tau_{ij}$ is the shear stress tensor, $q_j$ and $q_{v,j}$ are the $j$-component of the translational-rotational heat flux and vibrational-electronic heat flux, $V_{s,j}$ is the diffusion velocity of $s$th species,
\begin{equation}
  \tau_{il} = \mu(\frac{\partial u_i}{\partial x_l} + \frac{\partial u_l}{\partial x_i}) - \frac{2}{3}\mu \delta_{il} \frac{\partial u_k}{\partial x_k},
\end{equation}
\begin{equation}
  q_l = -k\frac{\partial T}{\partial x_l} + \rho \sum_{s=1}^{ns} h_s Y_s V_{s,l},
\end{equation}
\begin{equation}
  q_{v,l} = -k_v\frac{\partial T_v}{\partial x_l} + \rho \sum_{s=1}^{ns} e_{v,s} Y_s V_{s,l},
\end{equation}

\begin{equation}
  V_{s,l} = -D_s\frac{\partial Y_s}{\partial x_l},
\end{equation}
where $\mu, T, T_v$ are viscosity, translational-rotational temperature, and vibrational-rotational temperature. $k$ and $k_v,$ are the thermal conductivity of $T$ and $T_v$, respectively. $D_s$ is the effective diffusion coefficient of $s$th species. The total energy and vibrational enthalpy are given by,

\begin{equation}
  e_t = \frac{u_i u_i}{2} + \sum_{s=1}^{ns}h_s Y_s - \frac{p}{\rho},
\end{equation}
\begin{equation}
  e_v = \sum_{s=1}^{ns}e_{v,s} Y_s,
\end{equation}
where $h_s,e_{v,s}$ are the specific enthalpy and the specific vibrational enthalpy of $s$th species. The ideal gas assumption is introduced to close the system, 
\begin{equation}
  p = \rho \sum_{s=1}^{ns}Y_s \frac{R_u}{M_s} T,
\end{equation}
where $R_u$ is the universal gas constant, $M_s$ is the molecular weight of $s$th species. More detailed derivatives regarding the thermodynamic, transport, chemical reaction and vibrational excitation properties can be refer to Ref~\cite{SuNum,MARTINEZFERRER201488}.
%The choice of reaction mechanism of finite-rate chemistry depends on the composition of the ablation and pyrolysis products.  The governing equations are solved by our finite-volume in-house chemically reacting flow (CRF) solver~\cite{CAO2019679,SuNum}, which has been validated to be accurate and efficient in our previous work. The convective flux is splitted by the AUSMPW+ scheme~\cite{kim2001methods} and is reconstructed by the MUSCL scheme~\cite{VANLEER1979101} with the minmod limiter~\cite{SwebyHigh} to achieve a shock capturing capability. The 11 species ($\rm{N_2}$, $\rm{O_2}$, $\rm{NO}$, $\rm{C_3}$, $\rm{CO_2}$, $\rm{C_2}$, $\rm{CO}$, $\rm{CN}$, $\rm{N}$, $\rm{O}$, $\rm{C}$) with 18 reactions~\cite{doi:10.2514/1.J052659} is used for pure ablation employs  is taken from Ref.

\subsection{Implicit time integration}
The governing equations are written in the generalized coordinates, 
\begin{equation}
  \label{Eq.implicit}
  \frac{1}{J} \frac{\partial \bm{Q}}{\partial t} + \bm{R(Q)} = 0,
\end{equation}
where $J$ is the determinant of the curvilinear coordinate transformation Jacobian from Cartesian coordinate $(x1,x2,x3)$ to curvilinear coordinates $(\xi,\eta,\zeta)$. $\bm{R(Q)}$ is the residual associated with the spatial discretization,
\begin{equation}
  \label{Eq.iterate}
  \bm{R(Q)} = (\frac{\partial \bm{\tilde{E}}_{inv}}{\partial \xi} +  \frac{\partial \bm{\tilde{F}}_{inv}}{\partial \eta} +  \frac{\partial \bm{\tilde{H}}_{inv}}{\partial \zeta}) + \frac{Ma}{Re}(\frac{\partial \bm{\tilde{E}}_{vis}}{\partial \xi} +  \frac{\partial \bm{\tilde{F}}_{vis}}{\partial \eta} + \frac{\partial \bm{\tilde{H}}_{vis}}{\partial \zeta}) -\frac{1}{J}\cdot \bm{S},
\end{equation}
where $(\tilde{\bm{F}}_{inv}, \tilde{\bm{G}}_{inv},\tilde{\bm{H}}_{inv})$ and $(\tilde{\bm{F}}_{vis}, \tilde{\bm{G}}_{vis},\tilde{\bm{H}}_{vis})$ are the inviscid flues and viscid fluxes in the $(\xi,\eta,\zeta)$ direction, respectively. $Ma$ and $Re$ are the Mach number and Reynold number induced by nondimensionalization. The iteration of converging to a steady state (${\partial \bm{Q}}/{\partial t} = 0$) can be accelerated by employing local time-stepping $\Delta \tau$ in replace of physical time-stepping $\Delta t$ due to the fact that the choice of time step only impacts the characteristics of convergence, rather than the convergent solutions ($\bm{R(Q)} = 0$). For unsteady solutions, the dual time-stepping is usually employed for implicit scheme,%In this manuscript, the MUSCL scheme~\cite{VANLEER1979101} with the minmod limiter~\cite{SwebyHigh} is coupled with the Roe's upwind scheme to calculate the residual.
\begin{equation}
	\label{Eq.dualimplicit}
	\frac{1}{J}\frac{\partial \bm{Q}}{\partial\tau} + \frac{1}{J}\frac{\partial \bm{Q}}{\partial t} +\bm{R}(\bm{Q}) = 0
\end{equation}
where $\tau$ is the pseudo time. This dual time-stepping scheme is equal to Eq.\ref{Eq.implicit} when the pseudo iteration is convergent $\partial \bm{Q}/\partial \tau = 0$. The implicit iterate scheme of Eq.~\ref{Eq.dualimplicit} in differential form is written as, 
\begin{equation}
  (\frac{1}{J \tau} + \frac{1+ \theta}{J t}) \Delta \bm{Q}^{m+1} = \bm{R}(\bm{Q}^{m+1}) - \frac{(1+\theta)(\bm{Q}^{m+1}-\bm{Q}^n) - \theta(\bm{Q}^n - \bm{Q}^{n-1})}{J t},
\end{equation}
where $\Delta \bm{Q} = \bm{Q}^{m+1} - \bm{Q}^m$ is the increment of the conservative variables, supscript $m$ denotes the pseudo time level, $n$ denotes the physical time level. $\theta$ is the time differential coefficient that $\theta = 2$ for second-order time discretization and $\theta = 0$ for first-order time discretization. $\Delta \tau$ is computed by Courant-Friedrichs-Lewy (CFL) condition,
\begin{equation}
  \Delta \tau = \frac{J \cdot CFL}{(\lambda_{\xi,inv} + \lambda_{\eta,inv} + \lambda_{\zeta,inv}) + (\lambda_{\xi,vis} + \lambda_{\eta,vis} + \lambda_{\zeta,vis}) + \lambda_S},
\end{equation}
where $\lambda$ is the spectral radius (maximum eigenvalue) of flux jacobian, subscript $S$ denotes the source term. Linearizing $\bm{R}$ in Taylor expansion, 
\begin{equation}
  \label{Eq.linear}
  \bm{R(Q)}^{m+1} = \bm{R(Q)}^m + \frac{\partial \bm{R}}{\partial \bm{Q}}|^m \Delta \bm{Q}^{m+1} + O(\left\| \Delta \bm{Q}^{m+1}\right\|)^2,
\end{equation}
where $\partial \bm{R} / \partial \bm{Q}$ is the flux jacobian. Discarding the high-order truncation terms $O(\Delta \bm{Q}^{m+1})$ in Eq.(~\ref{Eq.linear}). Hence, the linearized implicit iteration scheme of Eq.(\ref{Eq.iterate}) is reformulated as 
\begin{equation}
  \label{Eq.linearizedimplicit}
  (\frac{1}{J \tau}\bm{I} + \frac{1+ \theta}{J t} \bm{I}+ \frac{\partial \bm{R}}{\partial \bm{Q}}) \Delta \bm{Q}^{m+1} = \bm{R}(\bm{Q}^{m}) - \frac{(1+\theta)(\bm{Q}^m-\bm{Q}^n) - \theta(\bm{Q}^n - \bm{Q}^{n-1})}{J t},
\end{equation}
where $\bm{I}$ is the identity matrix. Expanding Eq.(\ref{Eq.linearizedimplicit}) in generalized coordinates, 
\begin{equation}
	\label{Eq.implicitscheme}
  \begin{array}{l}
    \left\{   \bm{I}_m + J^{-1}\Delta \tau  [\partial_{\xi} (\bm{A}_{\xi,inv} + \bm{A}_{\xi,vis}) + \partial_{\eta} (\bm{A}_{\eta,inv} + \bm{A}_{\eta,vis}) + \partial_{\zeta} (\bm{A}_{\zeta,inv} + \bm{A}_{\zeta,vis}) ] \right\} \\ \Delta \bm{Q}^{m+1} = \bm{LHS} \Delta \bm{Q}^{m+1}  = \bm{RHS}^m
  \end{array}
\end{equation}
where $\partial_\xi, \partial_\eta, \partial_\zeta$ is the differential operator in the $\xi,\eta,\zeta$ direction. $\bm{LHS}$ and $\bm{RHS}$ refers to the left-hand side (implicit operator) and right-hand side, in relation to convergence characteristics and spatial residuals, respectively. $\bm{A}_{\xi,inv}, \bm{A}_{\eta,inv}$ and $\bm{A}_{\zeta,inv}$ are the flux jacobian of inviscid fluxes. $\bm{A}_{\xi,vis}, \bm{A}_{\eta,vis}$ and $\bm{A}_{\zeta,vis}$ are the flux jacobian of viscid fluxes. $\bm{A}_s$ is the flux jacobian of source term. These formulations are given by,
\begin{equation}
	\begin{array}{c}
		\bm{RHS}=\bm{R}-\frac{(1+\theta)(\bm{Q}^m-\bm{Q}^n) - \theta(\bm{Q}^n - \bm{Q}^{n-1})}{J t}, \bm{I}_m =  (\frac{1}{J\Delta \tau} + \frac{1+\theta}{J\Delta t}) \bm{I} - \bm{A}_s,\\
		\bm{A}_{\xi,inv} = \frac{\partial \tilde{\bm{F}}_{inv}}{\partial \bm{Q}},\bm{A}_{\eta,inv} = \frac{\partial \tilde{\bm{G}}_{inv}}{\partial \bm{Q}},\bm{A}_{\zeta,inv} = \frac{\partial \tilde{\bm{H}}_{inv}}{\partial \bm{Q}},\bm{A}_S = \frac{\partial \bm{S}}{J\partial \bm{Q}},\\
		\bm{A}_{\xi,vis} = \frac{Ma}{Re}\frac{\partial \tilde{\bm{F}}_{vis}}{\partial \bm{Q}},\bm{A}_{\eta,vis} = \frac{Ma}{Re}\frac{\partial \tilde{\bm{G}}_{vis}}{\partial \bm{Q}},\bm{A}_{\zeta,vis} = \frac{Ma}{Re}\frac{\partial \tilde{\bm{H}}_{vis}}{\partial \bm{Q}}.
		\end{array}
\end{equation}

The upwind spatial discretization schemes are commonly used to decompose the convective fluxes into positive and negative fluxes based on the eigenvalues of the flux Jacobian matrix, i.e., $\bm{\tilde{F}_{inv}} = \bm{A}\bm{Q} = \bm{S}^{-1} \bm{\Lambda} \bm{S} \bm{Q} = \bm{S}^{-1} (\bm{\Lambda}^++\bm{\Lambda}^-) \bm{S} \bm{Q}= \bm{A}^{+}\bm{Q}_{\mathrm{L}}+\bm{A}^-\bm{Q}_{\mathrm{R}}$, where $\mathrm{L,R}$ denote the left and the right variables, $\bm{\Lambda}$ and $\bm{S}$ are the eigenvalue and eigenvector matrix calculated by characteristic decomposition. Flux-difference splitting (FDS)~\cite{ROE1997250} and flux-vector splitting (FVS)~\cite{STEGER1981263} generates different approximations of flux jacobian. The first one is Roe splitting, $\bm{A}^{\pm} = (\bm{A} \pm | \tilde {\bm{A}} |)/2$, where $\tilde {\bm{A}}$ is the Roe averaged flux jacobian. The second one is eigenvalue splitting, $\bm{A}^{\pm} = (\bm{A} \pm \bm{S}^{-1} | \bm{\Lambda} | \bm{S})/2$. The third one is spectral splitting $\bm{A}^{\pm} = (\bm{A} \pm \lambda \bm{I})/2$, where $\lambda = \max{(\mathrm{diag}(\bm{\Lambda}))}$ is the spectral radius of the eigenvalue matrix. The spectral splitting is widely adopted in practical applications due to its high robustness and simple formats. Then, the differential form of the convective flux jacobian matrix in $\xi$ direction $\partial_{\xi} \bm{A}_{\xi,inv,i} = \bm{A}_{i+1/2} - \bm{A}_{i-1/2}$, resulting in
\begin{equation}
	\bm{A}_{\xi,inv,i-1} = -\bm{A}_{i-1/2}^+,\bm{A}_{\xi,inv,i} = \bm{A}_{i+1/2}^+ -\bm{A}_{i-1/2}^-,\bm{A}_{\xi,inv,i+1} = \bm{A}_{i+1/2}^-.
\end{equation}

The exact analytical form of viscous flux jacobian matrix is very complicated and is commonly approximated by viscous spectral,
\begin{equation}
  \begin{array}{c}
  \bm{A}_{\xi,vis,i-1}= -\lambda_{\xi,vis,i-1} \bm{I}, \bm{A}_{\xi,vis,i}= 2\lambda_{\xi,vis,i} \bm{I}, \bm{A}_{\xi,vis,i+1}= -\lambda_{vis,i+1} \bm{I}, \\ 
  \end{array}
  \label{Eq.vislambda}
\end{equation}

where $\lambda_{vis} = (\mu_l + \mu_t)\frac{Ma}{Re}\max({\frac{4\mu}{3\rho},\frac{k}{\rho c_v}},\max_{s=1}^{ns}(D_s)) \frac{\xi_x^2+\xi_y^2+\xi_z^2}{J\rho}$. Similar expressions hold for the $\eta$ and $\zeta$ direction. Expand Eq.\ref{Eq.implicitscheme} at point $(i,j,k)$,

\begin{equation}
  \label{Eq.ijkiterate}
  \begin{array}{c}
		(\bm{I}_{m} + \bm{A}_{\xi,i} + \bm{A}_{\eta,j}+ \bm{A}_{\zeta,k})\Delta \bm{Q}_{ijk} + \bm{A}_{\xi,i+1}\Delta \bm{Q}_{i+1} + \bm{A}_{\xi,i-1}\Delta \bm{Q}_{i-1} + \\ \bm{A}_{\eta,j+1}\Delta \bm{Q}_{j+1} +\bm{A}_{\eta,j-1}\Delta \bm{Q}_{j-1} + \bm{A}_{\zeta,k+1}\Delta \bm{Q}_{k+1} + \bm{A}_{\zeta,k-1}\Delta \bm{Q}_{k-1} =  \bm{RHS}_{ijk}^m
  \end{array}
\end{equation}

The above equation is a block sparse system and is solved by LUSGS method~\cite{yoon1988lower}. The LUSGS method decomposes the $\bm{LHS}$ in Eq.\ref{Eq.ijkiterate} into diagonal matrix $\bm{D}$, lower triangular matrix $\bm{L}$, and upper triangular matrix $\bm{U}$ before performing a factorization, 

\begin{equation}
	\begin{array}{c}
	(\bm{D} + \bm{L} + \bm{U}) \Delta \bm{Q} = (\bm{D}(\bm{I} + \bm{D}^{-1}\bm{L} + \bm{D}^{-1}))\Delta \bm{Q} \\ 
   \approx \bm{D}(\bm{I} + \bm{D}^{-1}\bm{L}) (\bm{I} +\bm{D}^{-1}\bm{U})\Delta\bm{Q}  = (\bm{D+L})\bm{D}^{-1}(\bm{D+U}) \Delta\bm{Q}
	\end{array}
\end{equation}

First, Eq.\ref{Eq.ijkiterate} is solved from front to back, 

\begin{equation}
  \begin{array}{c}
    (\bm{D} + \bm{L})\Delta \bm{Q}^{m+1,*} = \bm{RHS}^m, \\
    \bm{D} \bm{Q}^{m+1,*}  + \bm{A}_{i-1}\bm{Q}_{i-1}^{m+1,*} + \bm{A}_{j-1}\bm{Q}_{j-1}^{m+1,*}+ \bm{A}_{k-1}\bm{Q}_{k-1}^{m+1,*} =  \bm{RHS}^m,
  \end{array}
\end{equation}
where $\bm{Q}_{i-1}^{m+1,*},\bm{Q}_{j-1}^{m+1,*}$, and $\bm{Q}_{k-1}^{m+1,*}$ have been solved in the previous iteration and can be moved to the right-hand side. Second, Eq.\ref{Eq.ijkiterate} is solved from back to front,
\begin{equation}
  \begin{array}{c}
    (\bm{D} + \bm{U})\Delta \bm{Q}^{m+1} = \bm{D}\Delta\bm{Q}^{m+1,*}, \\
    \bm{D} \bm{Q}^{m+1}  + \bm{A}_{i+1}\bm{Q}_{i+1}^{m+1} + \bm{A}_{j+1}\bm{Q}_{j+1}^{m+1}+ \bm{A}_{k+1}\bm{Q}_{k+1}^{m+1} =  \bm{D}\Delta\bm{Q}^{m+1,*},
  \end{array}
\end{equation}
where $\bm{Q}_{i+1}^{m+1},\bm{Q}_{j+1}^{m+1}$, and $\bm{Q}_{k+1}^{m+1}$ have been solved in the previous iteration and can be moved to the right-hand side.% The increments $\Delta \bm{Q}^{m+1}$ are used to update $\bm{Q}$ until $\bm{RHS}$ decreases to a specified magnitude.

\subsection{Component-splitting method}

The above implicit time integration uses an unified spectral radius in flux splitting, proclaimed as the coupled implicit method. However, density/momentum/energy equations and component equations have different characteristic wave speeds, i.e., different eigenvalue properties, which could impact the convergence speed. The component-splitting method separates the implicit operator into flow part solving density/momentum/energy equations and component part solving component equations. Differential flux jacobian is reformulated corresponding to flow part and component part for accelerated convergence. Different spectral radius of viscid fluxes are employed corresponding to flow equations and component equations, rather than using a unified viscous spectral. The convective flux jacobian is separated into flow part and component part, discarding the cross terms.

\begin{equation}
	\label{Eq.convectiveJacobian}
	\begin{aligned}
		\bm{A}_{\xi,inv} = \frac{1}{J}\left[{\begin{array}{*{20}{c}}
    0\\ 
    \beta E_k \xi_x - Uu \\
    \beta E_k \xi_y - Uv \\
    \beta E_k \xi_z - Uw \\
    \left[(\gamma-2)E_k  - h\right]U \\
    -Ue_v\\
    -UY_1\\
    \ldots\\
    -UY_{ns-1}
    \end{array}}\quad
    {\begin{array}{*{20}{c}}
    \xi_x\\
    A_x\\
    B_y\\
    C_z\\
    \epsilon_x\\
    \xi_xe_v\\
    Y_1\xi_x\\
    \ldots\\
    Y_{ns-1}\xi_x
    \end{array}}\quad
    {\begin{array}{*{20}{c}}
    \xi_y\\
    B_x\\
    A_y\\
    B_z\\
    \epsilon_y\\
    \xi_ye_v\\
    Y_1\xi_y\\
    \ldots\\
    Y_{ns-1}\xi_y
    \end{array}}\quad
    {\begin{array}{*{20}{c}}
    \xi_z\\
    C_x\\
    C_y\\
    A_z\\
    \epsilon_z\\
    \xi_ze_v\\
    Y_1\xi_z\\
    \ldots\\
    Y_{ns-1}\xi_z
    \end{array}}
    {\begin{array}{*{20}{c}}
    0\\
    \beta \xi_x\\
    \beta \xi_y\\
    \beta \xi_z\\
    \gamma U\\
    0\\
    0\\
    \ldots\\
    0
    \end{array}} \right. & \\
		\left. \quad{\begin{array}{*{20}{c}}
      0\\
      \Phi \xi_x\\
      \Phi \xi_y\\
      \Phi \xi_z\\
      \Phi U\\
      U\\
      0\\
      \ldots\\
      0
      \end{array}}\quad{\begin{array}{*{20}{c}}
      0\\
      \tilde{\xi}_1\xi_x\\
      \tilde{\xi}_1\xi_y\\
      \tilde{\xi}_1\xi_z\\
      \tilde{\xi}_1 U\\
      -Ue_v\\
      U\\
      \ldots\\
      0
      \end{array}}\quad{\begin{array}{*{20}{c}}
    \ldots\\
    \ldots\\
    \ldots\\
    \ldots\\
    \ldots\\
    \ldots\\
    \ldots\\
    \ldots\\
    \ldots
    \end{array}}\quad \begin{array}{*{20}{c}}
    0\\
    \tilde{\xi}_{ns-1}\xi_x\\
    \tilde{\xi}_{ns-1}\xi_y\\
    \tilde{\xi}_{ns-1}\xi_z\\
    \tilde{\xi}_{ns-1} U\\
    -Ue_v\\
    0\\
    \ldots\\
    U
    \end{array}  \right]&
	\end{aligned}
\end{equation}

where $U = u\xi_x + v\xi_y + w\xi_z + \xi_t, \beta = \gamma -1, \Phi = -\beta + \rho_{electron}/M_{electron},  E_k = (u^2+v^2+w^2)/2, \tilde{\xi_s} = \gamma RT/M_s-\beta h_s, A_x = U+(2-\gamma)u\xi_x,B_x = -\beta v\xi_x+u\xi_y,A_y = U+(2-\gamma)v\xi_y,B_y = -\beta u\xi_y+v\xi_x,A_z = U+(2-\gamma)w\xi_z,B_z = -\beta v\xi_z+w\xi_y, \epsilon_x = (E_k + h)\xi_x - \beta U u, C_x = -\beta w \xi_x+u \xi_z, \epsilon_y = (E_k + h)\xi_y -\beta U v, C_y = -\beta w\xi_y+v\xi_z, \epsilon_z = (E_k + h)\xi_z - \beta U w,C_z = -\beta u\xi_z+w \xi_x$. Characteristic decomposition on this flux jacobian matrix generates the eigenvalue matrix $\bm{\Lambda} = diag(U-c\nabla\xi,U,U,U,U+c\nabla\xi,U,\cdots,U)/J$ and the spectral radius $\lambda_{inv} = |U| + c\nabla \xi$, where $\nabla \xi=\sqrt{\xi_x^2+\xi_y^2+\xi_z^2}$.

The implicit operator of component equations and flow equations can be independently solved without specific order. Discarding the cross effect between flow equations and component equations in implicit operator introduces component-splitting errors during iterations. Thus, all $ns$ component equations are solved with consistence correction to address numerical inconsistencies of mass fraction.
% \subsubsection{\label{sec:level3}Third-level heading: Citations and Footnotes}
\subsubsection{Implicit operator for flow equations} 

The implicit operator of flow equations solves the conservative variables $\bm{Q}_F=[\rho,\rho u, \rho v, \rho w, \rho e]^{\mathrm{T}}$ and freeze species mass fraction. The interacted flux jacobian between flow equations and component equations are neglected. Thus, the flux jacobian $\bm{A}$ between convective flux and specie density is 0 and remains the upper left $[6\times 6]$ block of Eq.\ref{Eq.convectiveJacobian}.
\begin{equation}
  \bf{A} = \left[{\begin{array}{*{20}{c}}
    0\\ 
    \beta E_k \xi_x - Uu \\
    \beta E_k \xi_y - Uv \\
    \beta E_k \xi_z - Uw \\
    \left[(\gamma-2)E_k  - h\right]U \\
    -Ue_v\\
    \end{array}}\quad
    {\begin{array}{*{20}{c}}
    \xi_x\\
    A_x\\
    B_y\\
    C_z\\
    \epsilon_x\\
    \xi_xe_v\\
    \end{array}}\quad
    {\begin{array}{*{20}{c}}
    \xi_y\\
    B_x\\
    A_y\\
    B_z\\
    \epsilon_y\\
    \xi_ye_v\\
    \end{array}}\quad
    {\begin{array}{*{20}{c}}
    \xi_z\\
    C_x\\
    C_y\\
    A_z\\
    \epsilon_z\\
    \xi_ze_v\\
    \end{array}}
    {\begin{array}{*{20}{c}}
    0\\
    \beta \xi_x\\
    \beta \xi_y\\
    \beta \xi_z\\
    \gamma U\\
    0\\
    \end{array}}\quad{\begin{array}{*{20}{c}}
      0\\
      \Phi \xi_x\\
      \Phi \xi_y\\
      \Phi \xi_z\\
      \Phi U\\
      U\\
      \end{array}}
    \right] 
\end{equation}

The spectral splitting ($\bm{A}^{\pm} = \bm{A} \pm \lambda \bm{I}$) is used to generate the differential form of the flux jacobian matrix. The spectral radius of is evaluated by $\lambda = \lambda_{inv} + \lambda_{vis}$, $\lambda_{inv} = (|U| + c \nabla \xi)/J$, $\lambda_{vis} = (\mu_l + \mu_t)\frac{Ma}{Re}\max({\frac{4\mu}{3\rho},\frac{k}{\rho c_v}}) \frac{\nabla \xi^2}{J\rho}$. Then, the differential form of the flux jacobian is
\begin{equation}
	\label{Eq.fluxjacobiandiff}
\bm{A}_{\xi,i} = (\lambda_{inv,i} + 2\lambda_{vis,i})\bm{I}, \bm{A}_{\xi,i-1} = -\bm{A}_{i-1/2}^+, \bm{A}_{\xi,i+1} = \bm{A}_{i+1/2}^-.
\end{equation}

\subsubsection{Implicit operator for component equations} 

The implicit operator of component equations solves only the conservative variables $\bm{Q}_C=[\rho Y_1, \rho Y_2, \cdots, \rho Y_{ns}]^{\mathrm{T}}$ keeping $\bm{Q}_F$ unchanged. Thus the flux jacobian $\bm{A}$ remains the lower right $[ns-1 \times ns-1]$ block of Eq.\ref{Eq.convectiveJacobian} and extends it to $ns$ species, i.e., $U \bm{I}_{ns \times ns}$. The spectral radius of component equations is given by $\lambda = \lambda_{inv} + \lambda_{vis}$, $\lambda_{inv} = (|U|)/J$, $\lambda_{vis} = (\mu_l + \mu_t)\frac{Ma}{Re}\max_{s=1}^{ns}(D_s) \frac{\nabla \xi^2}{J\rho}$. The differential form of flux jacobian is given by
\begin{equation}
	\label{Eq.fluxjacobiandiff}
	\bm{A}_{\xi,i} = (\lambda_{inv,i} + 2\lambda_{vis,i})\bm{I}, \bm{A}_{\xi,i-1} = -(U_{i-1}+\lambda_{i})\bm{I}, \bm{A}_{\xi,i+1} = (U_{i+1} -\lambda_{i})\bm{I}.
\end{equation}

Using the spectral splitting in flux splitting, the flux Jacobian matrix has all diagonal elements equal, while the off-diagonal elements are zero. In consequence, Eq.\ref{Eq.ijkiterate} can be solved in scalar operation instead of matrix operation to avoid exponential increases in the computational consumption with the number of species. This makes the time cost of implicit time integration acceptable for large number of species in multicomponent equations. However, the computational cost for detailed chemistry and characteristic decomposition with large number of species remains unaffordable. A remedy for this issue can be partial decomposition~\cite{WANG2019364} and tabulated chemistry~\cite{VICQUELIN20111481}, which can be coupled with component-splitting method to solve multicomponent reactive flows. 

\subsubsection{Consistence correction}

The coupled implicit scheme solves only $ns-1$ component equations and uses $\rho_{ns} = \rho - \sum_{s=1}^{ns-1}\rho_{s}$ to avoid numerical inconsistencies of mass fraction. However, this correction is not suitable for component splitting method, as the errors accumulate in the last species. To address this issue, all $ns$ specie equations are solved by introducing two types of consistence correction to avoid numerical inconsistencies that the total mass fraction not equaling 1. The first correction normalizes the increment of conservative variables,
\begin{equation}
	\rho_{s}^{m+1} = \rho_{s}^{m} +\Delta Q_{\rho_s} + W_s(\Delta Q_{\rho} - \sum_{s=1}^{ns}\Delta Q_{\rho_s}),
\end{equation}
where $W_s$ is the weight coefficient, $\Delta Q_{\rho_s}, \Delta Q_{\rho}$ are the increment of specie density and total density, respectively. The weight coefficient is computed by $W_s = Q_{\rho_s}/\sum_{s=1}^{ns}\Delta {Q}_{\rho_s}$.Obviously, summing up the above equation yields $\rho^{m+1} = \rho^{m} +\Delta Q_{\rho}$, which is equivalent to the flow part. This correction only modifies the increment without altering the conservative variables at current time level, thereby ensuring conservation properties. 

The second correction normalizes the specie mass fractions,
\begin{equation}
	\rho_{s}^{m+1} = (\rho^{m} + \Delta Q_{\rho})\frac{\rho_s^m + \Delta Q_{\rho_s}}{\sum_{s=1}^{ns}(\rho_s^m +\Delta Q_{\rho_s})}.
\end{equation} 

Using consistence correction finally obtains the increment of conservative variables $\Delta \bm{Q}$. Since accuracy is determined by the grid resolution, spatial discretization method, and residual magnitude, this consistence correction does not compromise the accuracy when the residual achieves convergence. In fact, it can potentially improve the accuracy by reducing the residual to lower magnitude. %If the conservative variables appear with negative values, it will be assigned a value of 0 directly to continue iterations in both the coupled implicit method and the component-splitting method.

\section{Numerical cases}
% 前两个算例（激波管和对流扩散）验证精度，后面算例验证加速效果
Hypersonic flows are commonly in the state of thermo-chemical non-equilibrium due to the real gas effects~\cite{anderson2000hypersonic} and solve multicomponent Navier-Stokes equations for numerical simulations\cite{10.1063/5.0045184}. Numerical tests on hypersonic flows are conducted to assess the performance of the component-splitting method in accelerating implicit iteration. The efficiency, stability and computation time are compared between the component-splitting (CS) scheme with the coupled implicit (CI) scheme, where CS-1 and CS-2 use the first and the second consistence correction, respectively. All numerical cases employ LUSGS solving the block-sparse system and are performed on a dual-socket EPYC 7742 desktop.

\subsection{Time acceleration test}
\begin{figure}[htb!]
	\centering
		\includegraphics[width=0.45\textwidth]{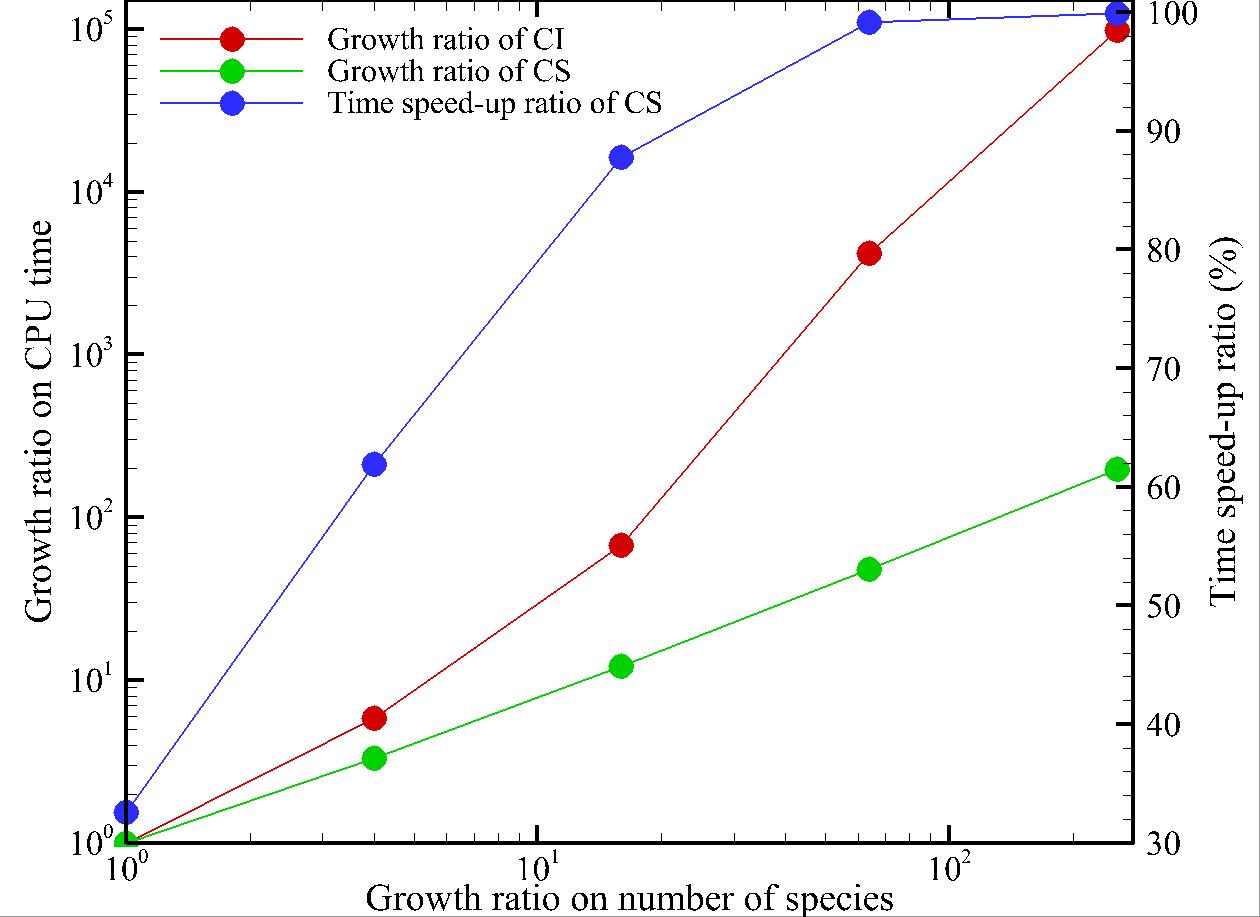}
	\caption{Growth ratio on the CPU time and speed-up ratio with respect to the growth ratio on the number of species from $ns=16$ to $ns=4096$.}
	\label{FIG.timeaceleration}
\end{figure}

The time acceleration ratio on implicit time integration of the CS method with respect to the number of species number is investigated through non-reactive uniform flows in a cubic box with a grid size of $10 \times 10 \times 10$, and all boundaries set to far-field equaling to the uniform flows. The chemical reaction is not considered to exclude the time cost of the detailed chemistry. The CPU time consumption of implicit time integration in a single iteration step, including assembling of flux jacobian matrix and solving the block-sparse system, is counted with the number of species increased from 16 to 4096. Let $t_{\mathrm{CS}}(ns)$ being the CPU time consumption of implicit time integration using the CS method to solve the multicomponent equations with $ns$ species, the growth ratio on CPU time $(t(ns)/t(16))$ and the speed-up ratio ($t_{\mathrm{CI}}(ns)-t_{\mathrm{CS}}(ns)$)$/t_{\mathrm{CI}}(ns)$ with respect to the growth ratio on the number of species $(ns/16)$ are show in FIG.~\ref{FIG.timeaceleration}. The computation time of the CS method is 33\% and 1\% of that of the CI method for 16 and 1024 species, respectively, demonstrating an acceleration in computational efficiency with increasing number of species. The time consumption of CS increases linearly with the number of species, whereas the CI method exhibits exponential growth. As a result, the CS method enhances computational efficiency and ensures the feasibility of the implicit time integration for multicomponent reactive equations when using a large number of species.
% \begin{table}[htb!]
% 	\centering
% 	\caption{Comparison of CPU time (s) of single iteration between CS and CI method.}
% 	\label{Tab.CPUtime-ns}%添加标题 设置标签
% 	\begin{tabular}{cccc}
% 		\toprule
% 		%\hline
% 		species number&CPU time (s) of CI&CPU time (s) of CS&speed-up ratio (\%)\\
% 		\midrule
% 		% 16&$2.8\times 10 ^{-5}$&$2.5\times 10 ^{-5}$&
% 		% 64&$7.7\times 10 ^{-5}$&$3.9\times 10 ^{-5}$&
% 		% 256&$8.11\times 10 ^{-4}$&$1.51\times 10 ^{-4}$&
% 		% 1024&$3.3816\times 10 ^{-2}$&$3.41\times 10 ^{-4}$&
% 		% 4096&$7.10351\times 10 ^{-1}$&$1.24\times 10 ^{-3}$&
% 		% 16384&
% 		16&$6.694\times 10^{-3}$&$4.5100\times 10^{-3}$&32.6\\
% 		64&$3.9227\times 10^{-2}$&$1.4959\times 10^{-2}$&61.9\\
% 		256&$4.4936\times 10^{-1}$&$5.4932\times 10^{-2}$&87.8\\
% 		1024&$2.8190\times 10^{1}$&$2.1695\times 10^{-1}$&99.2\\
% 		4096&$6.5762\times 10^{2}$&$8.9015\times 10^{-1}$&99.9\\
% 		%\hline
% 		\bottomrule
% 	\end{tabular}
% \end{table}

\subsection{Steady hypersonic flows on cylinder}
\begin{figure}[htb!]
	\centering
	\subfigure
	{
		\includegraphics[height=0.3\textwidth]{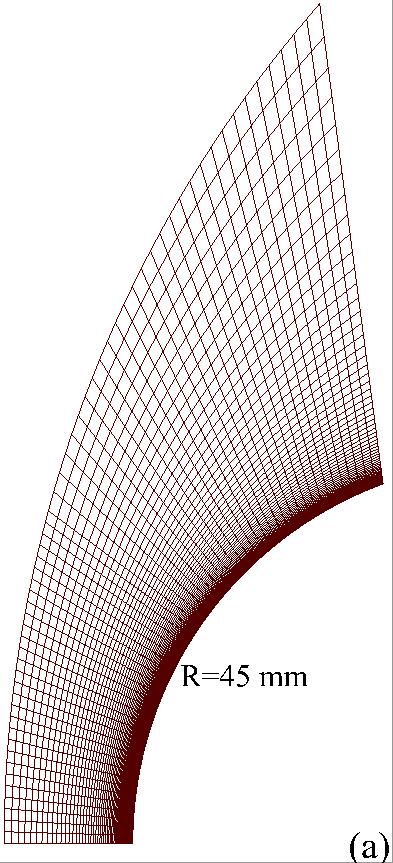}
	}
	\subfigure
	{
		\includegraphics[height=0.3\textwidth]{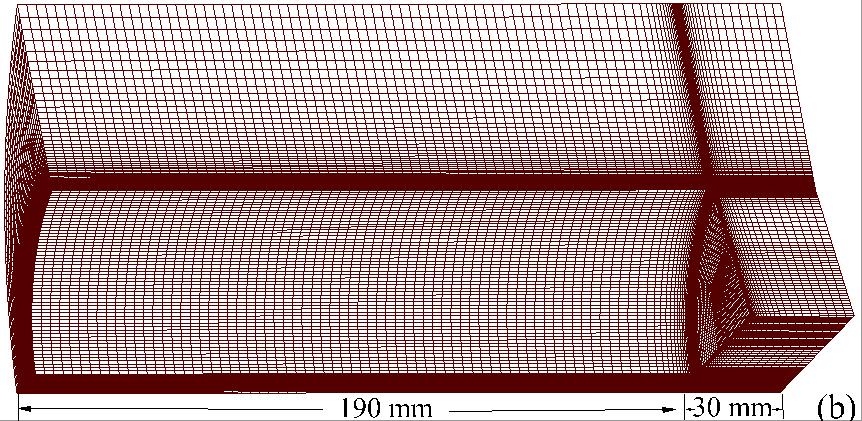}
	}
	\centering
	\caption{Computational grid of cylinder with radius of 45 mm:(a)cylinder2D; (b)cylinder3D.}
	\label{FIG.cylindergrid}
\end{figure}

High Enthalpy Shock Tunnel Gottingen (HEG)~\cite{doi:10.2514/6.2003-4252} conducted serval experiments on a 90mm diameter cylinder to measure the density distribution in the shock layer as well as surface pressure and heat flux. These experimental data have become a basis validation for numerical methods of thermo-chemical nonequilibrium flows. The free stream conditions are obtained from experimental inflow conditions and are presented in Table.~\ref{Tab.cylindercondition}. Numerical simulations are performed in both 2D and 3D, where the 2D case utilized a mesh grid of $60\times 95$ points in the circumferential and normal directions, respectively. The 3D case is conducted to simulate the experimental setup in wind tunnel with spanwise extension of 190 mm and a distance of 30 mm from the boundary. The spanwise direction is distributed with 110 grid points. The 2D and 3D computational mesh are shown in the FIG.~\ref{FIG.cylindergrid}. The wall boundary is assumed to be non-slip, non-catalytic, and isothermal at 300 K. The Park model~\cite{park1989assessment} which contains 11 species with 21 elementary reactions is used to simulate chemical reactions of the air mixture with the reaction rate constants are obtained from Ref.\cite{gnoffo1989conservation}. The MUSCL scheme~\cite{VANLEER1979101} coupled with the Roe scheme\cite{ARABI2019178} is applied to calculate the convective flux, where the Minod limiter~\cite{SwebyHigh} is used to achieve a shock capturing capability. The Spalart-Allmaras (S-A) turbulence model~\cite{doi:10.2514/6.1992-439} and second order central difference scheme are employed to calculate the viscous fluxes.

\begin{table}[htb!]
	\centering
	\caption{Free-stream conditions for numerical simulations of cylinder.}
	\label{Tab.cylindercondition}%添加标题 设置标签
	\begin{tabular}{cccccccc}
		\toprule
		%\hline
		$Ma$&$p$(Pa)&$T(\mathrm{K})$&$Y_{N2}$&$Y_{O2}$&$Y_{N}$&$Y_{O}$&$Y_{NO}$\\
		\midrule
		8.98&476&901&0.7543&0.00713&$6.5 \times 10 ^{-7}$&0.2283&0.01026\\
		%\hline
		\bottomrule
	\end{tabular}
\end{table}

\begin{figure}[htb!]
	\centering
	\subfigure
	{
		\includegraphics[width=0.45\textwidth]{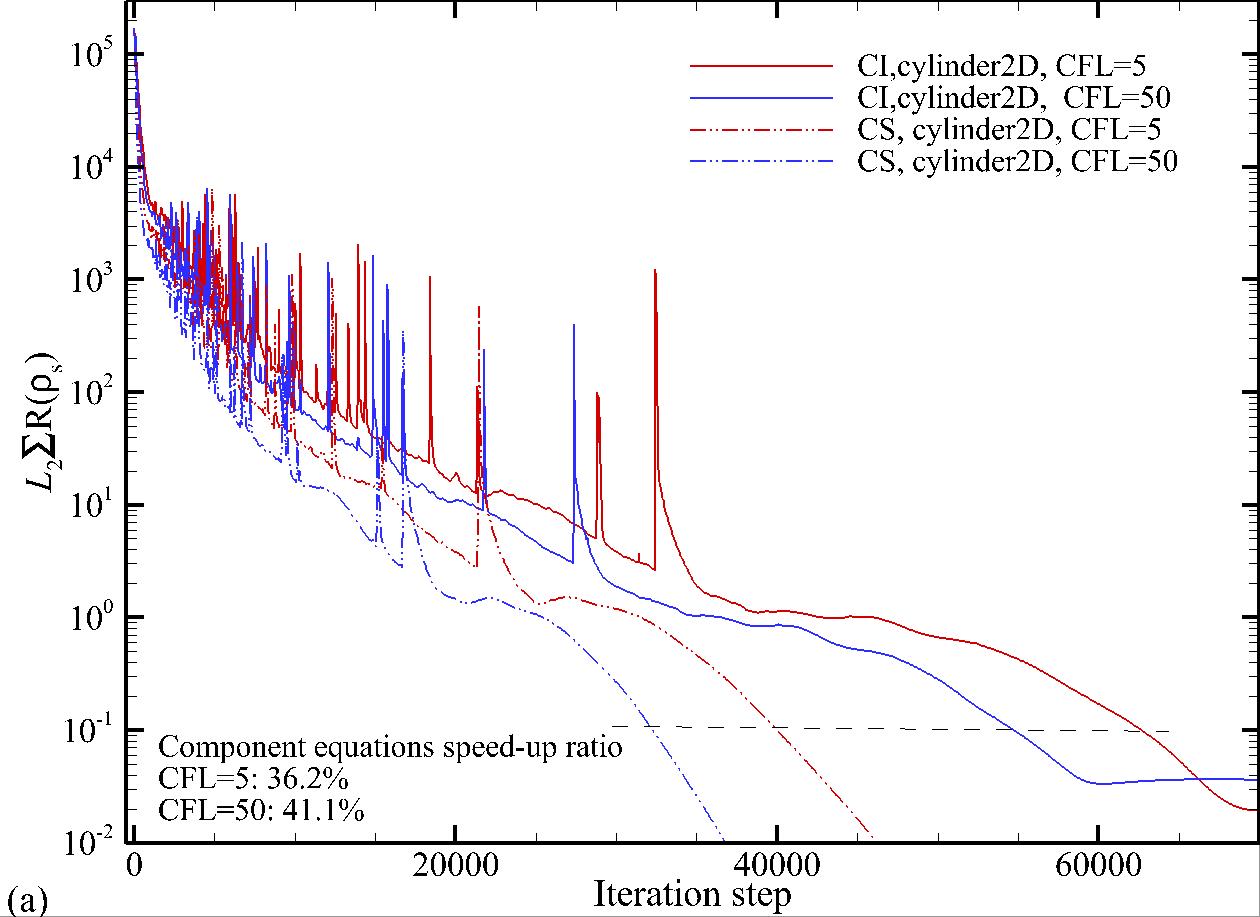}
	}
	\subfigure
	{
		\includegraphics[width=0.45\textwidth]{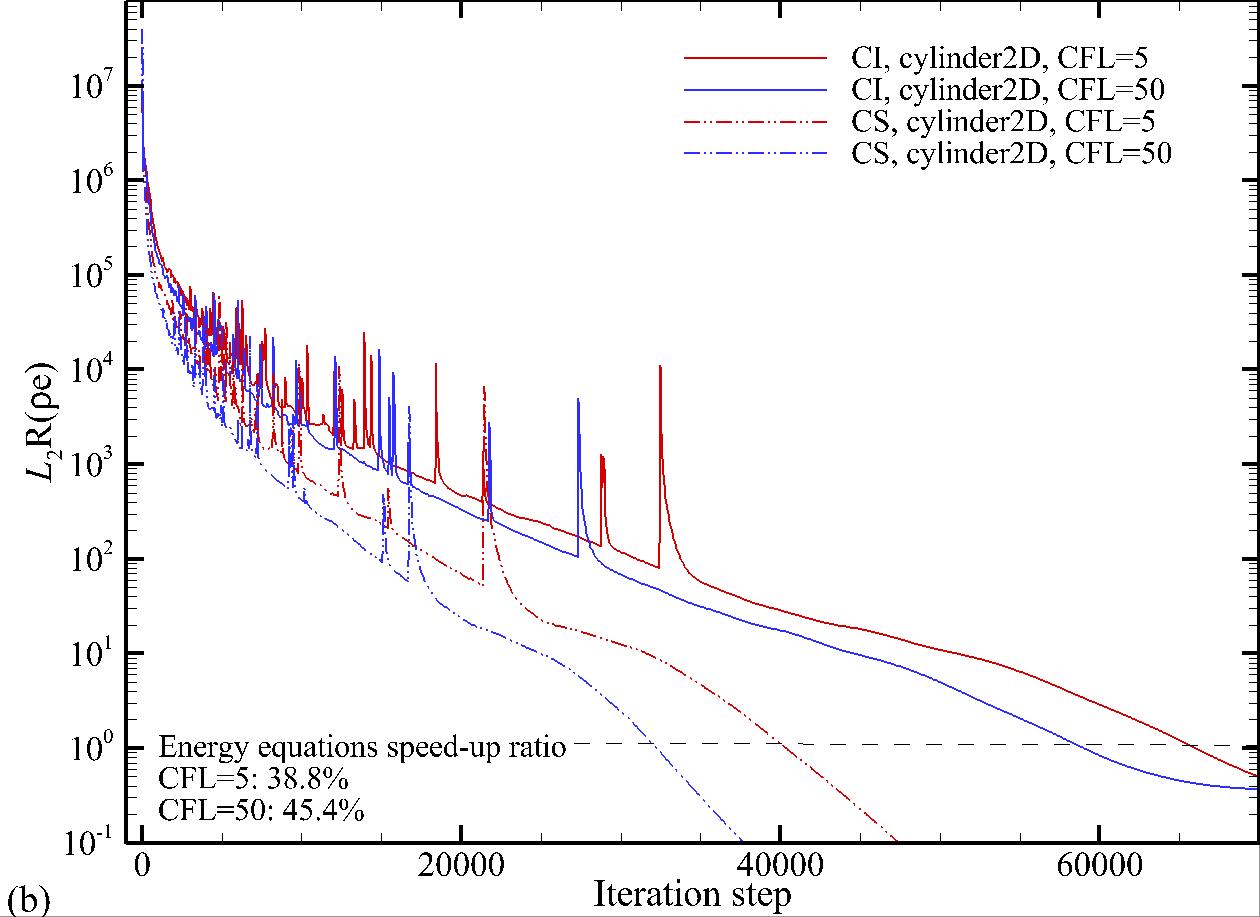}
	}
	\subfigure
	{
		\includegraphics[width=0.45\textwidth]{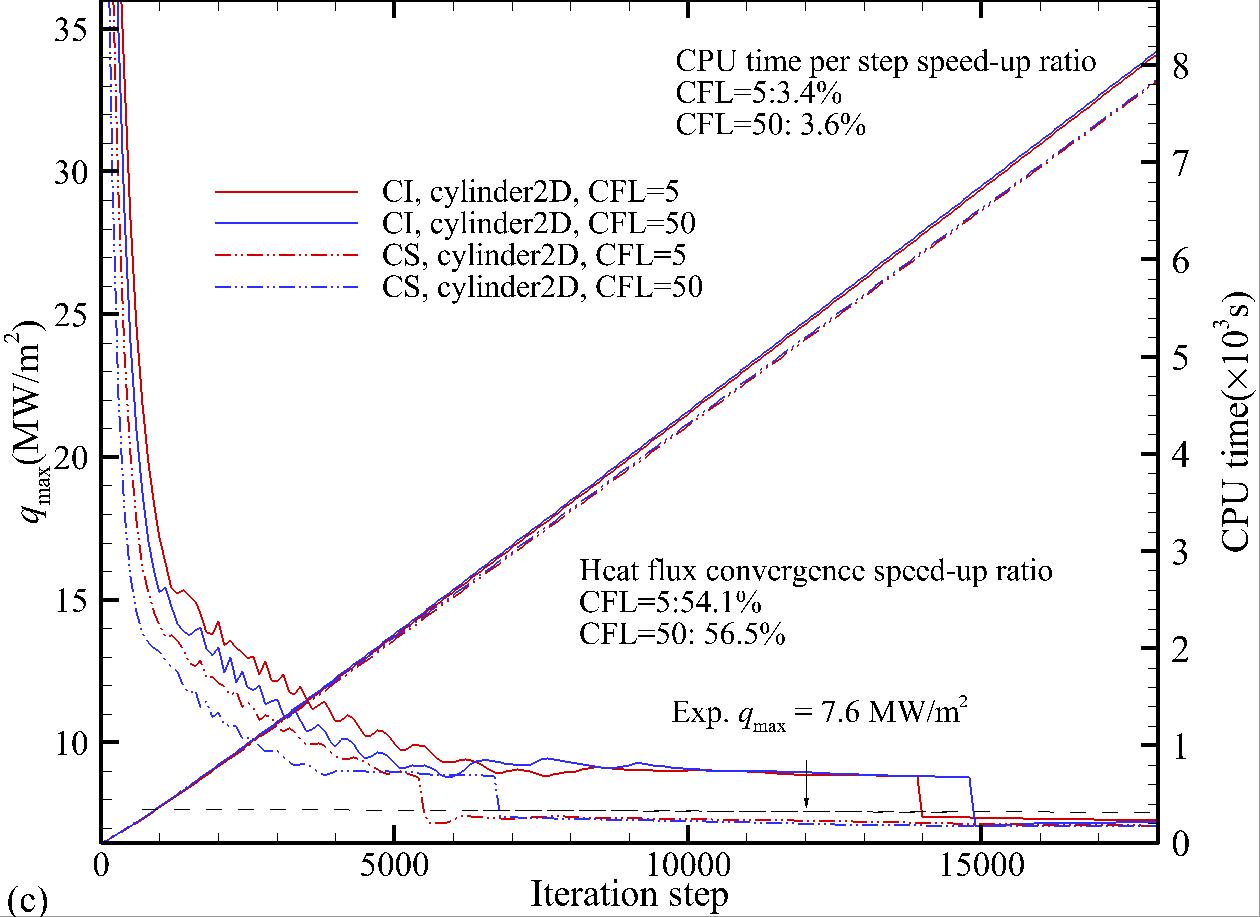}
	}
	\subfigure
	{
		\includegraphics[width=0.45\textwidth]{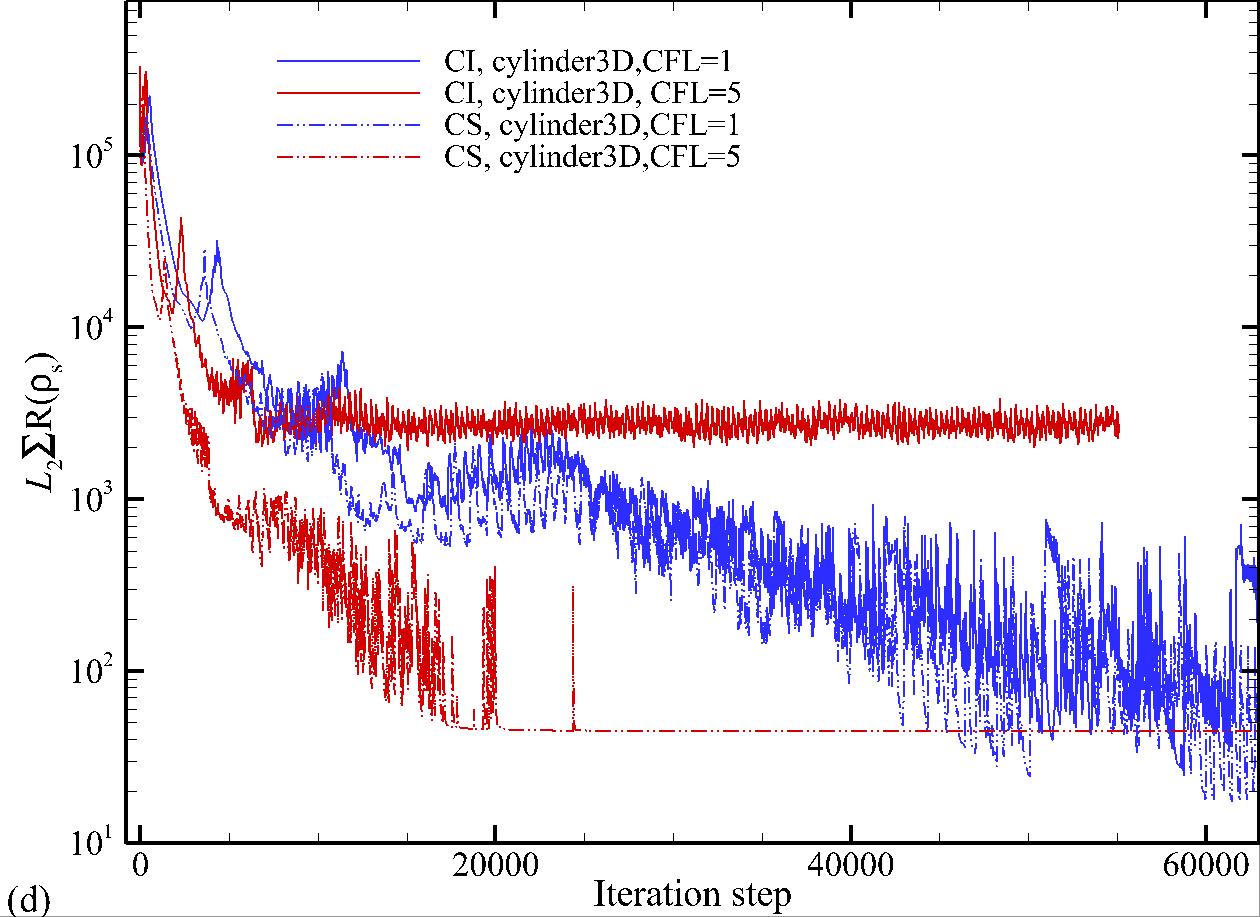}
	}
	\subfigure
	{
		\includegraphics[width=0.45\textwidth]{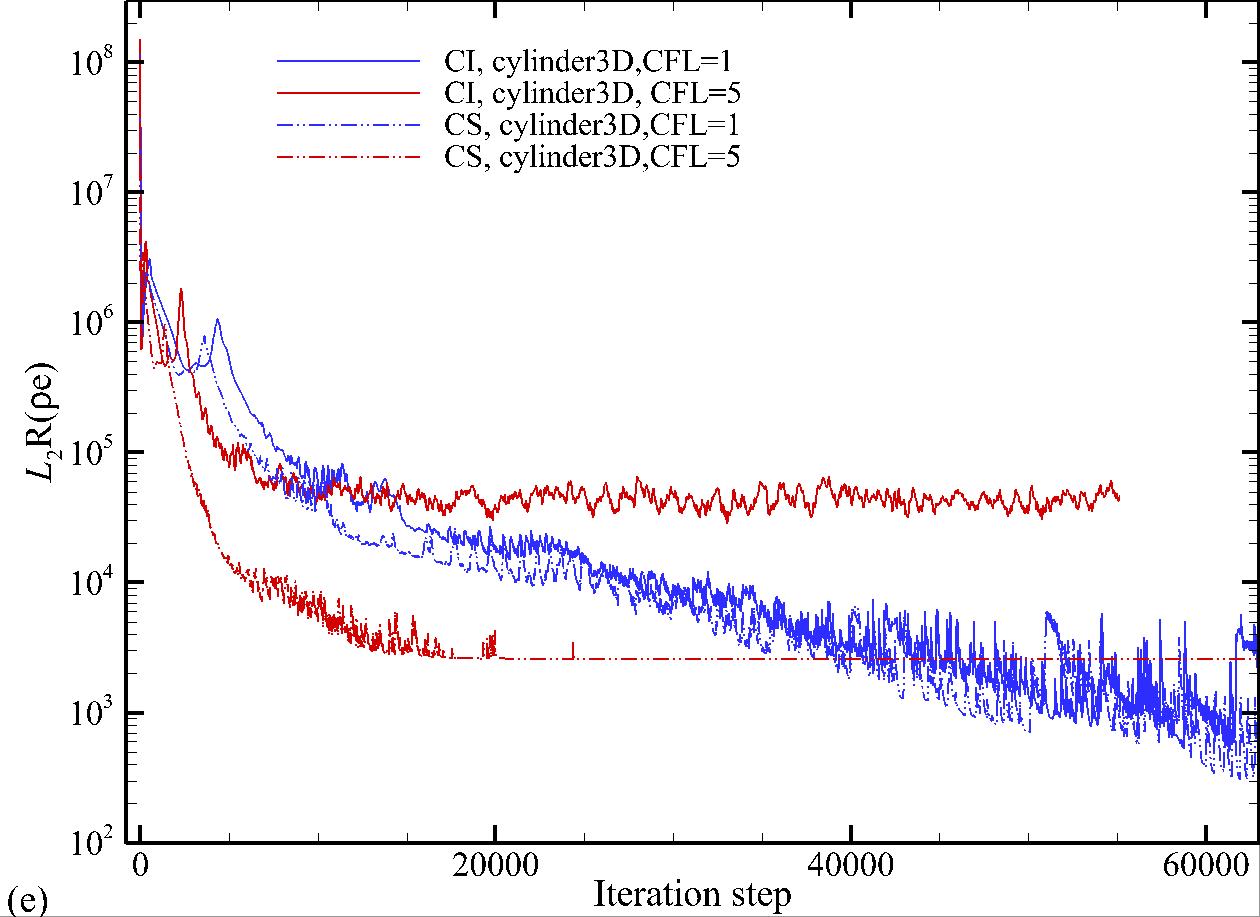}
	}
	\subfigure
	{
		\includegraphics[width=0.45\textwidth]{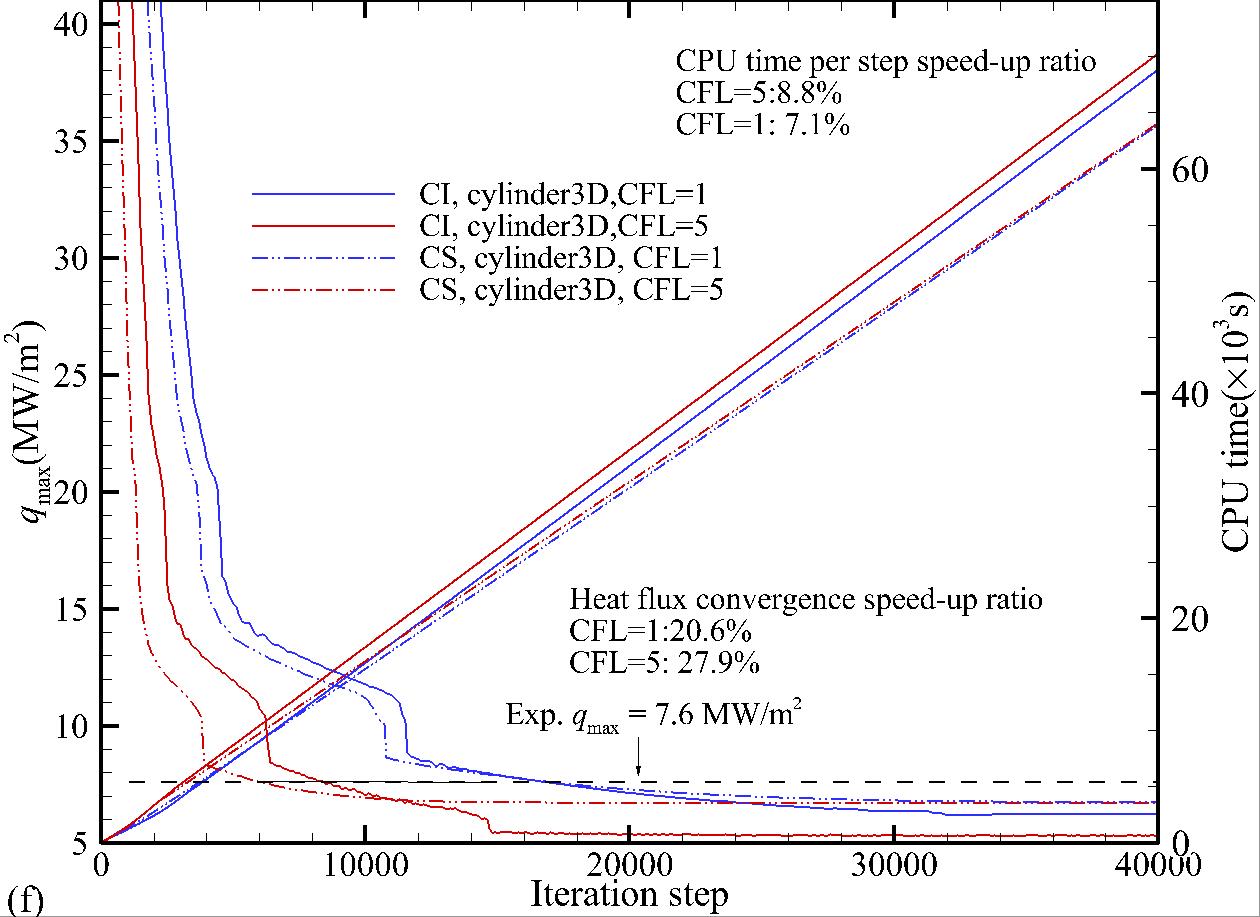}
	}
	\caption{Convergence curve of residual and heat flux with respect to iteration step for comparison of the CS and CI method in the cylinder case: ((a)(d)) residual of specie density; ((b)(e)) residual of energy; ((c)(f)) wall heat flux at stagnation points and CPU time cost; ((a)(b)(c)) cylinder2D; ((d)(e)(f)) cylinder3D. The dashed lines represent convergence criteria.}%$L_2$ norm of the sum of the residual of the specie density among all grid cells
	\label{FIG.cylinderconvergence}
\end{figure}

\begin{figure}[htb!]
	\centering
	\subfigure
	{
		\includegraphics[width=0.45\textwidth]{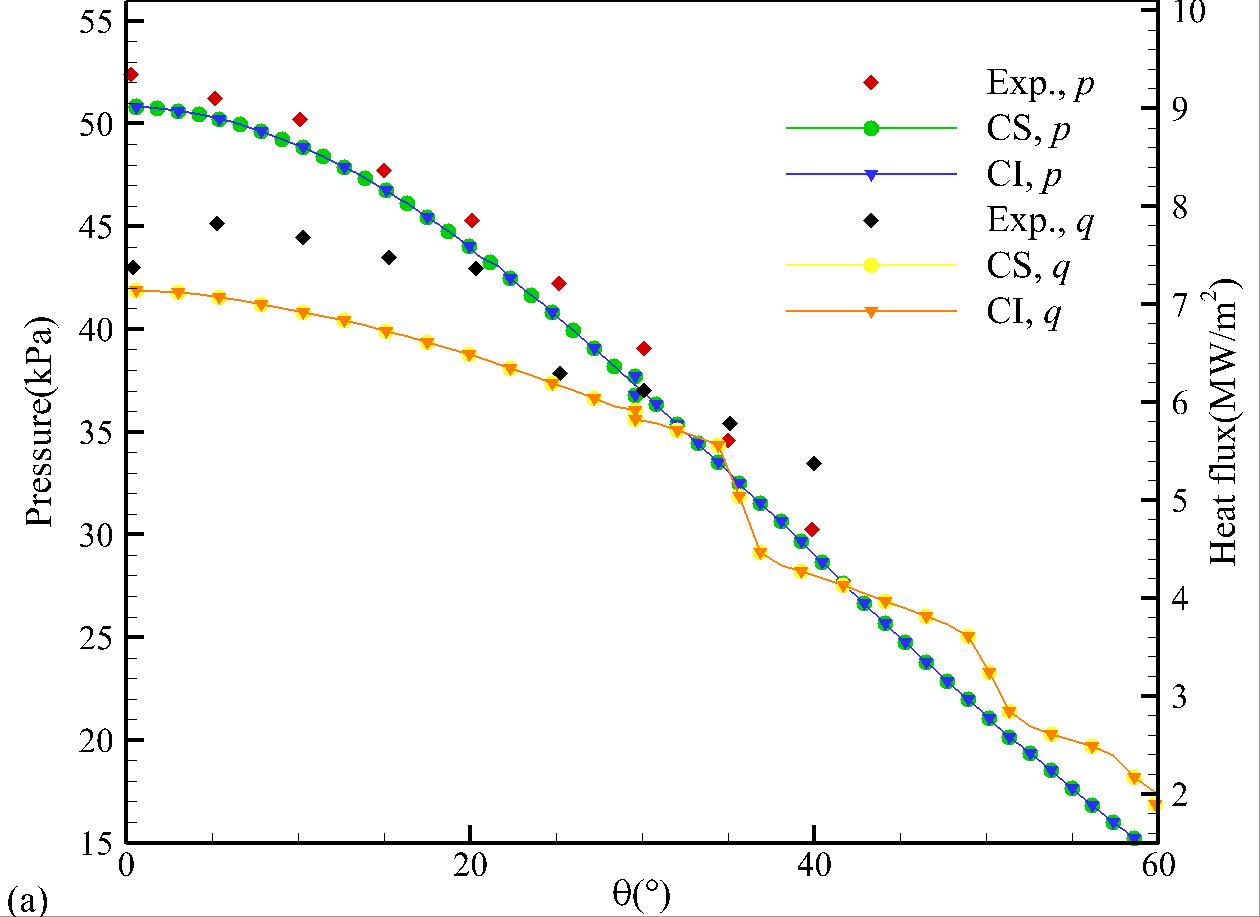}
	}
	\subfigure
	{
		\includegraphics[width=0.45\textwidth]{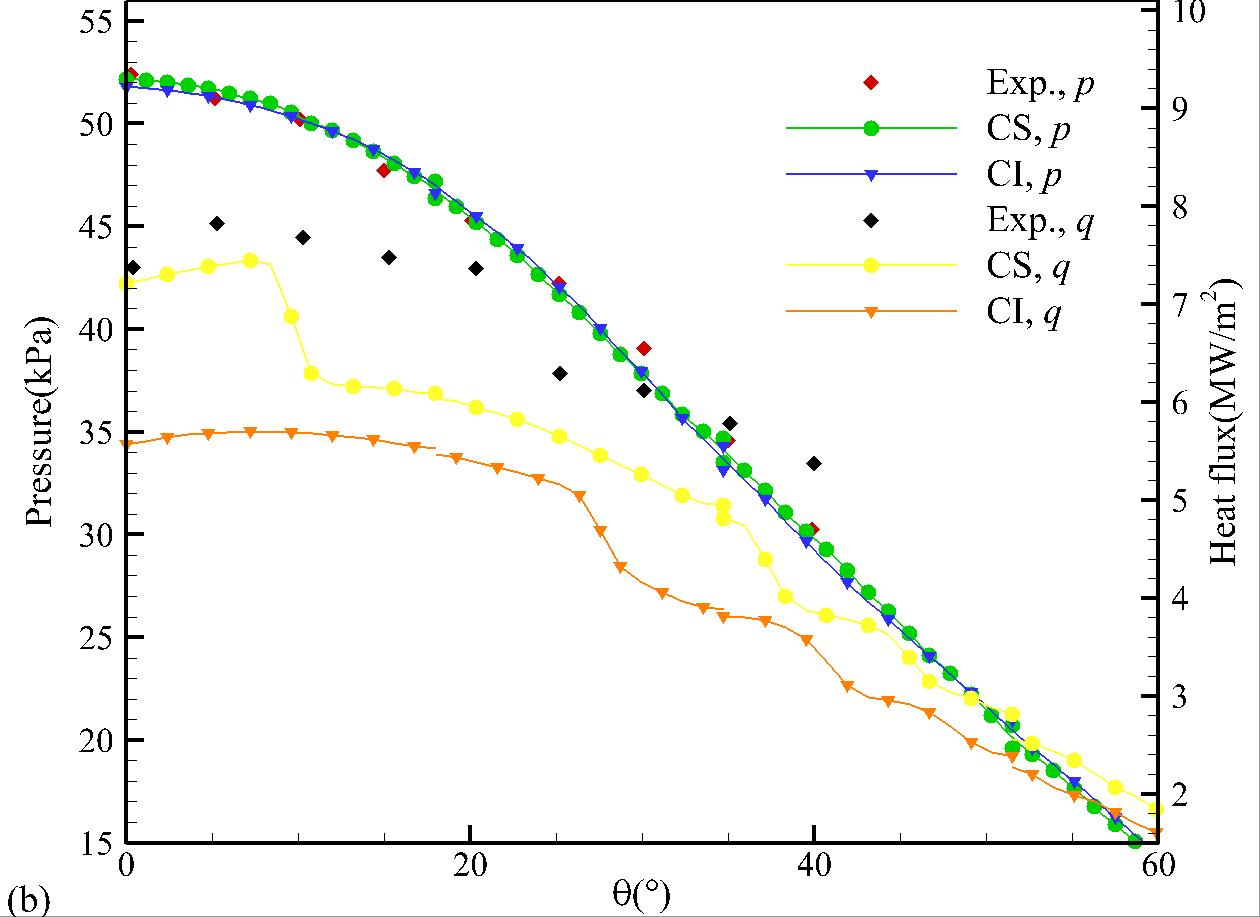}
	}
	\caption{Comparison of wall pressure and heat flux between CS, CI, and experimental data~\cite{KNIGHT20128} for hypersonic flows on cylinder.(a):cylinder2D;(b)cylinder3D}
	\label{FIG.cylinder-p-q}
\end{figure}

The convergence characteristics of residual and heat flux are compared between component-splitting method and coupled implicit method. The convergence curves of residual and wall heat flux are provided in FIG.~\ref{FIG.cylinderconvergence}. In 2D simulations, the convergence criterion is set to residual magnitude reaching $10^{-1}$, the CS method approximately accelerates the convergence of species density by 36.2\% and 41.1\% for CFL=5 and CFL=50, respectively. The convergence characteristics of the energy equation are similar to those of the component equations with the speed-up ratio is roughly equal. After considering the time consumption for flux calculations, input/output, etc., the CPU time per iteration step is reduced by 3.4\%-3.6\% using the CS method. With a convergence criterion of a relative error less than 1\% for heat flux, the required iteration step to convergence is reduced by 54.1\% and 56.5\% for CFL=5 and CFL=50, respectively. In 3D simulations, the CS method can use higher CFL number for better convergence of residual of both species density and energy while diminishing the residual to lower magnitude. The computed wall heat flux at stagnation points of CS and CI, after convergence, are in close agreement with the experimental results~\cite{doi:10.2514/6.2003-4252}, with a relative error of approximately 6\%. FIG.\ref{FIG.cylinder-p-q} compares the wall pressure and heat flux between the results from CS, CI, and experimental data~\cite{KNIGHT20128,10.1063/5.0047341}. In 2D simulations, due to the convergence of residual obtained from CS and CI method, errors are primarily attributed to spatial discretization. The results of the CS and CI methods are essentially identical due to the convergence of residual, indicating that the CS method does not compromise accuracy. In 3D simulations, the CS method achieves higher accuracy in wall heat flux calculations due to its lower magnitude of residual, indicating superior convergence characteristics of the component-splitting method. 

\subsection{Axisymmetric Shock-Wave Boundary Layer Interaction(ASWBLI)}
NASA Turbulence Modeling Resources (TMR)~\cite{doi:10.2514/6.2015-0316} present a basis validation for hypersonic flows. This validation is based on experiments~\cite{Kussoy1989Documentation} on a cone/ogive cylinder with a radius of $10.15$ cm a flare angle of 20$^{\circ}$. The inflow conditions are summarized in Table\ref{Tab.ASBLIcondition}. The computational domain is a quarter model with $161\times 201\times 31$ grid points in streamwise, transverse, and circumferential direction. The wall boundary assumes non-slip, non-catalytic, and isothermal of 311K. The spatial discretization method and chemical kinetic mechanism used to calculate the residual remains consistent with the approach employed in the cylinder case. To compare the efficiency of the component-splitting method with different numbers of species, $ns=5$ non-ionized model ($\rm{N_2}$, $\rm{O_2}$, $\rm{NO}$, $\rm{N}$, $\rm{O}$) and $ns=11$ ionized model ($\rm{N_2}$, $\rm{O_2}$, $\rm{NO}$, $\rm{N}$, $\rm{O}$, $\rm{NO^{+}}$, $\rm{N_2^{+}}$, $\rm{O_2^{+}}$, $\rm{N^+}$, $\rm{O^+}$, $\rm{e^-}$) are used for computation. 

\begin{table}[htb!]
	\centering
	\caption{Free-stream conditions for numerical simulations of ASBLI case.}
	\label{Tab.ASBLIcondition}%添加标题 设置标签
	\begin{tabular}{cccccc}
		\toprule
		%\hline
		Re&$Ma$&AOA($^{\circ}$)&$T(\mathrm{K})$&$Y_{N2}$&$Y_{O2}$\\
		\midrule
		$5.706 \times 10^{6}$&7.11&0.0&80&0.767&0.233\\
		%\hline
		\bottomrule
	\end{tabular}
\end{table}

\begin{figure}[htb!]
	\centering
	\subfigure
	{
		\includegraphics[width=0.45\textwidth]{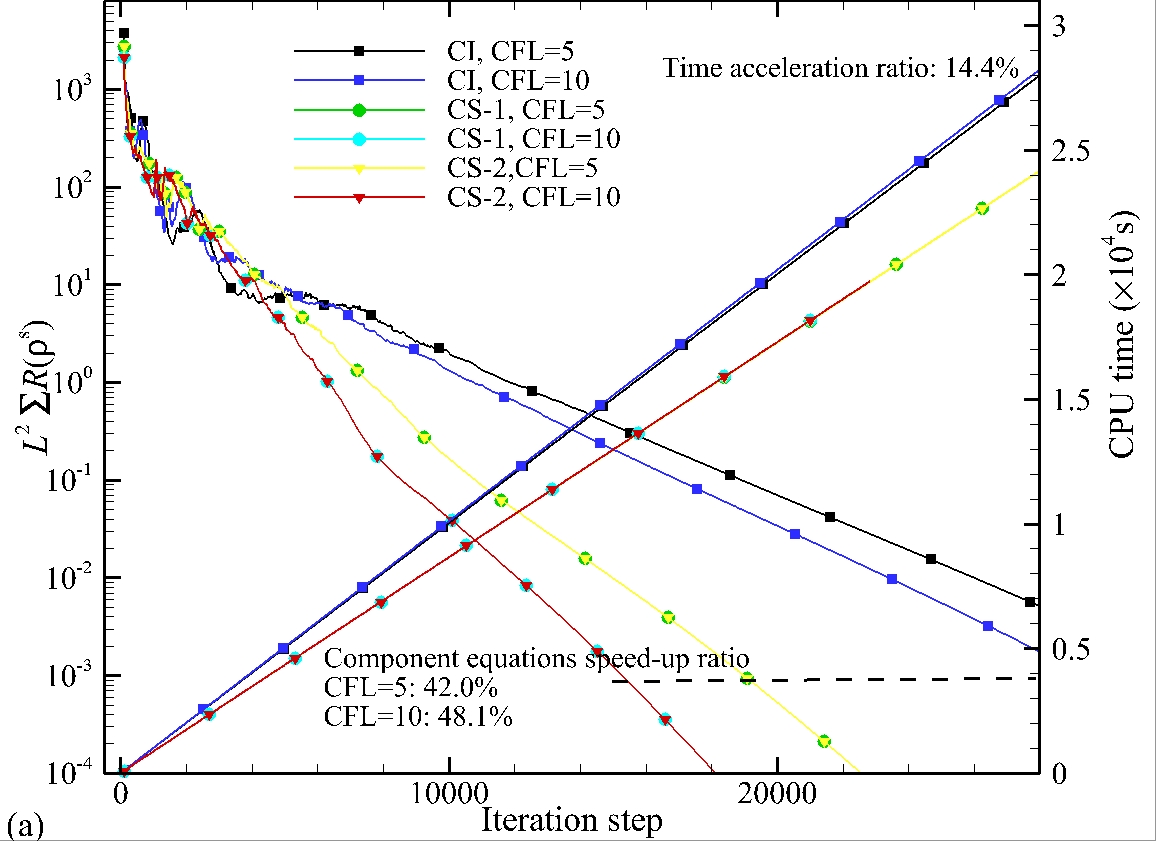}
		% \label{FIG.GSC}
	}
	\subfigure
	{
		\includegraphics[width=0.45\textwidth]{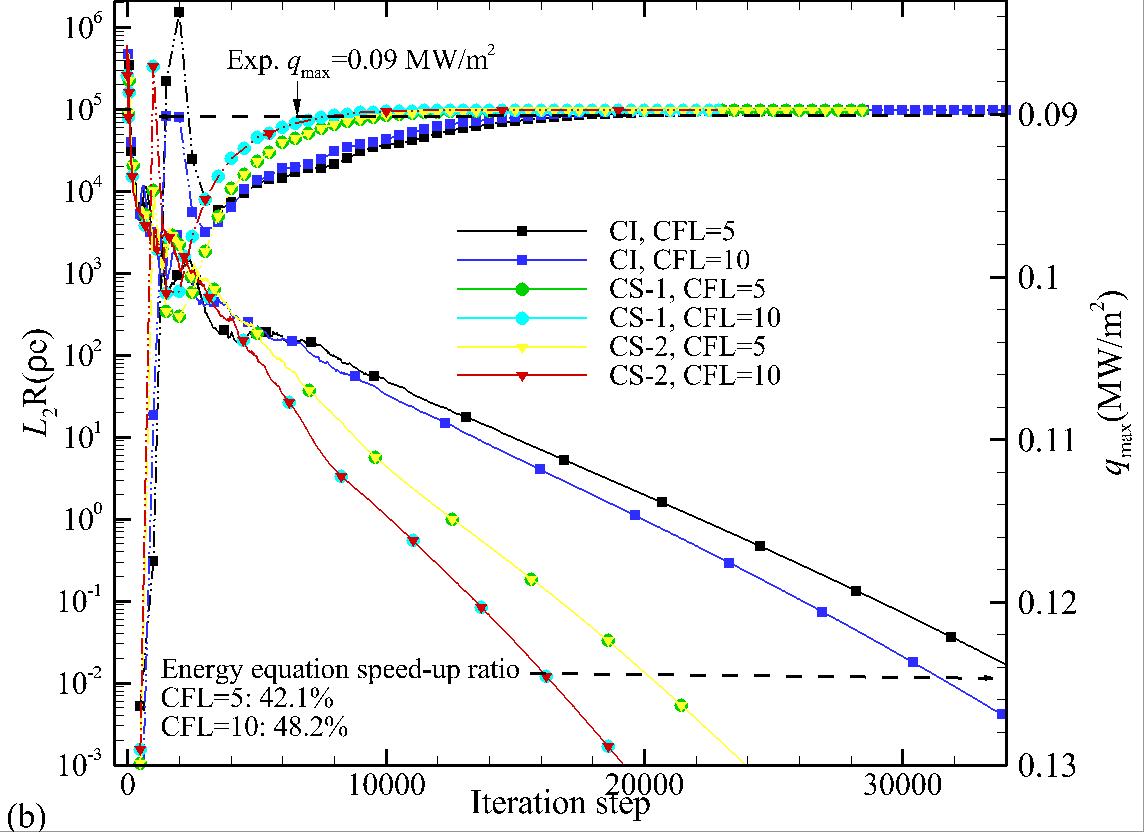}
		% \label{FIG.GSC_grid_a}
	}
	\subfigure
	{
		\includegraphics[width=0.45\textwidth]{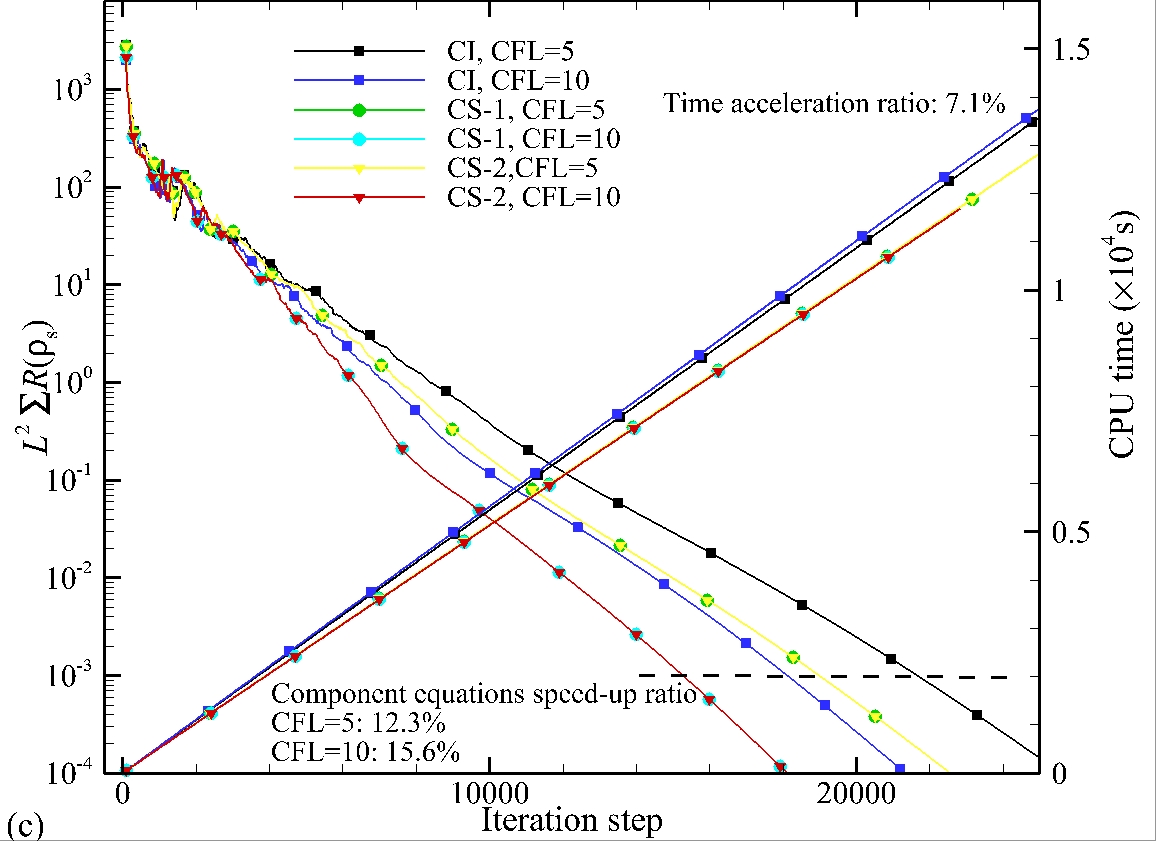}
		% \label{FIG.GSC}
	}
	\subfigure
	{
		\includegraphics[width=0.45\textwidth]{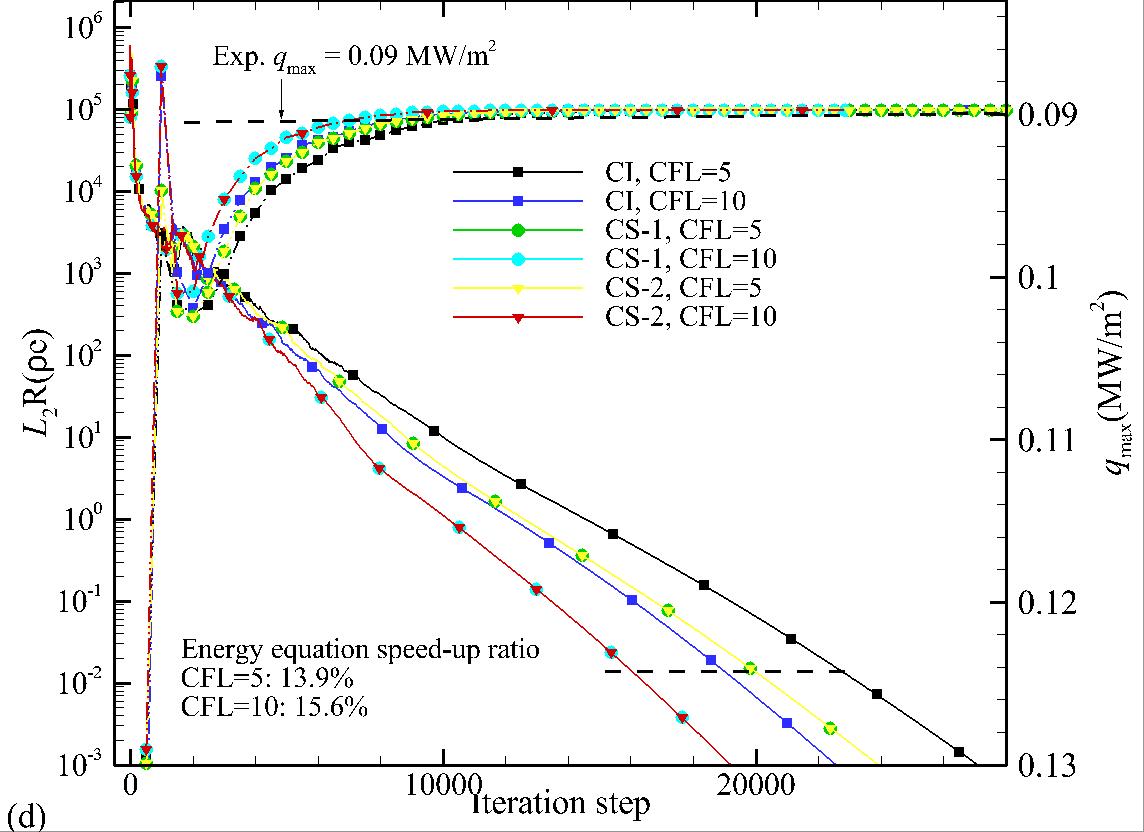}
		% \label{FIG.GSC_grid_a}
	}
	\caption{Convergence curve of the residual:((a),(c))$L_2$ norm of sum of residual of specie density and CPU time; ((b),(d))$L_2$ norm of sum of residual of energy and wall heat flux at stagnation points;((a),(b))$ns=11$ case;((c),(d))$ns=5$ case.}
	\label{FIG.ASWBLIcurve}
\end{figure}

The convergence curve of residual of species density and energy are shown in FIG.\ref{FIG.ASWBLIcurve}. The residual convergence behavior of the component equations is similar to that of the residual of the energy equation. In $ns=11$ case, the CS method achieve speed-up ratio of 42.0\% and 48.1\% for CFL=5 and CFL=10, respectively. Larger CFL number can improve the convergence acceleration effect of component-splitting method. With the convergence criterion set as relative error less than 1\%, the convergence of heat flux at stagnation points is accelerated by 42.8\% and 48.0\% for CFL=5 and CFL=10, respectively. Comparing the convergence curve between $ns=11$ case and $ns=5$ case, the convergence acceleration effects of residual and heat flux are enhanced with the increase number of species. The CPU time of single step is reduced by 14.4\% and 7.1\% for $ns=11$ case and $ns=5$ case, respectively.

\begin{figure}[htb!]
	\centering
	\subfigure
	{
		\includegraphics[width=0.3\textwidth]{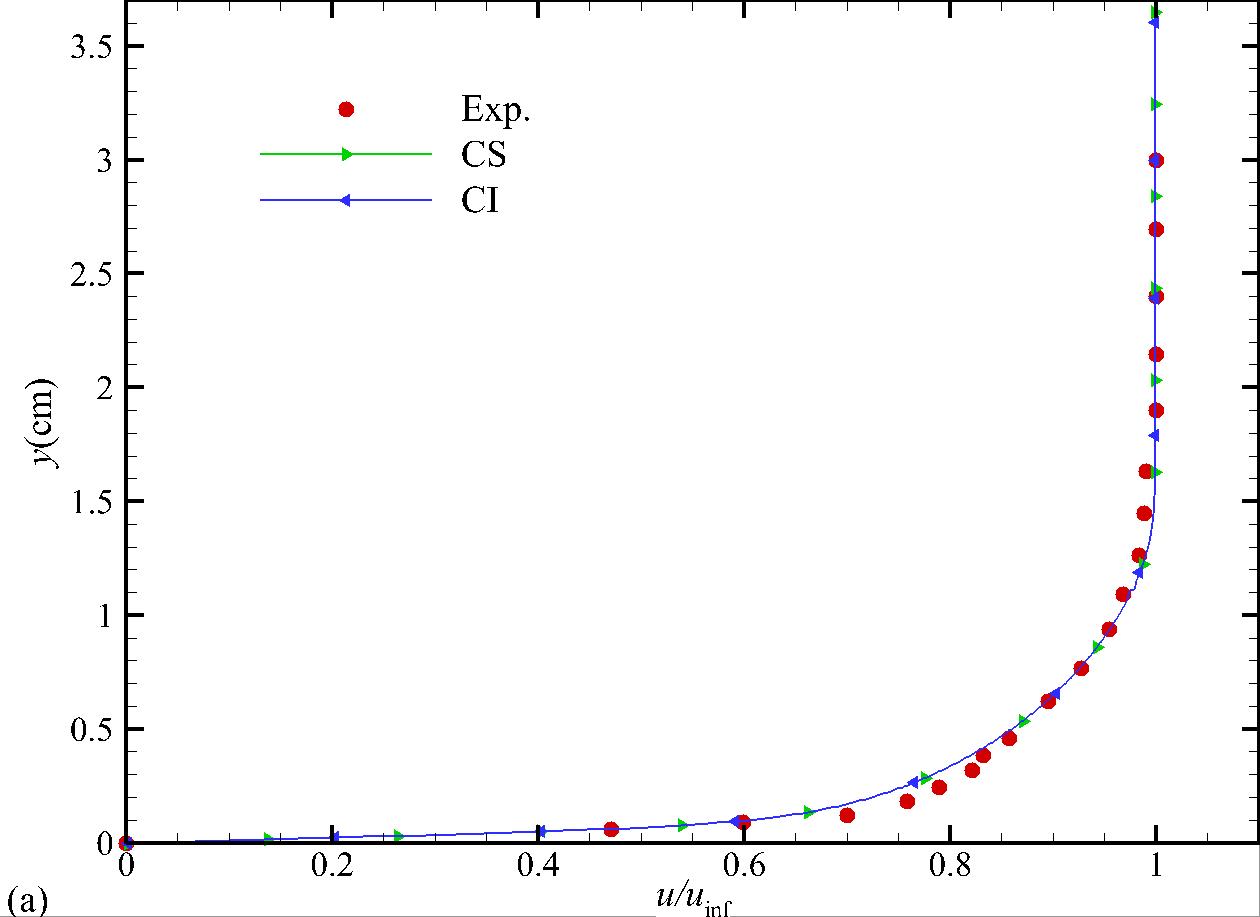}
		% \label{FIG.GSC}
	}
	\subfigure
	{
		\includegraphics[width=0.3\textwidth]{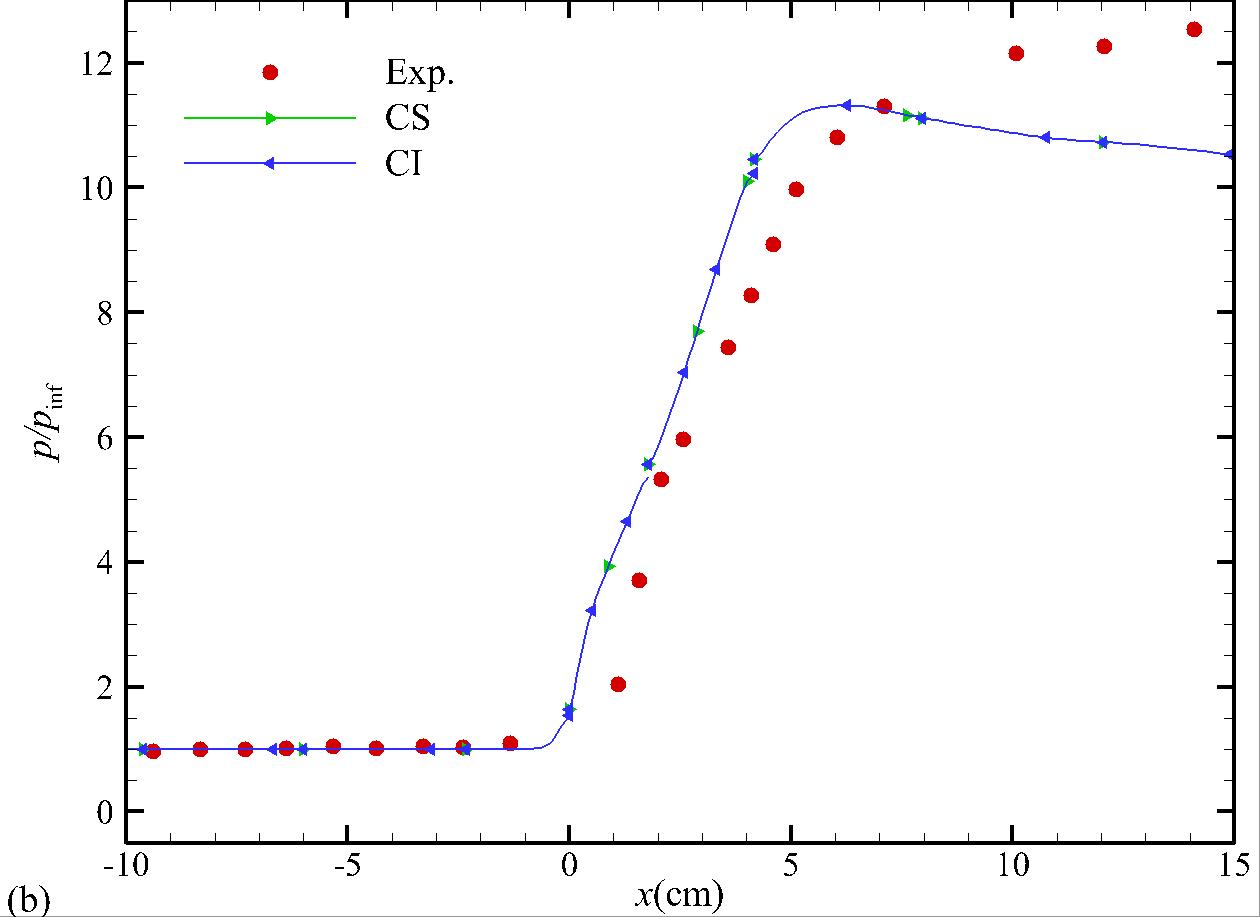}
		% \label{FIG.GSC_grid_a}
	}
	\subfigure
	{
		\includegraphics[width=0.3\textwidth]{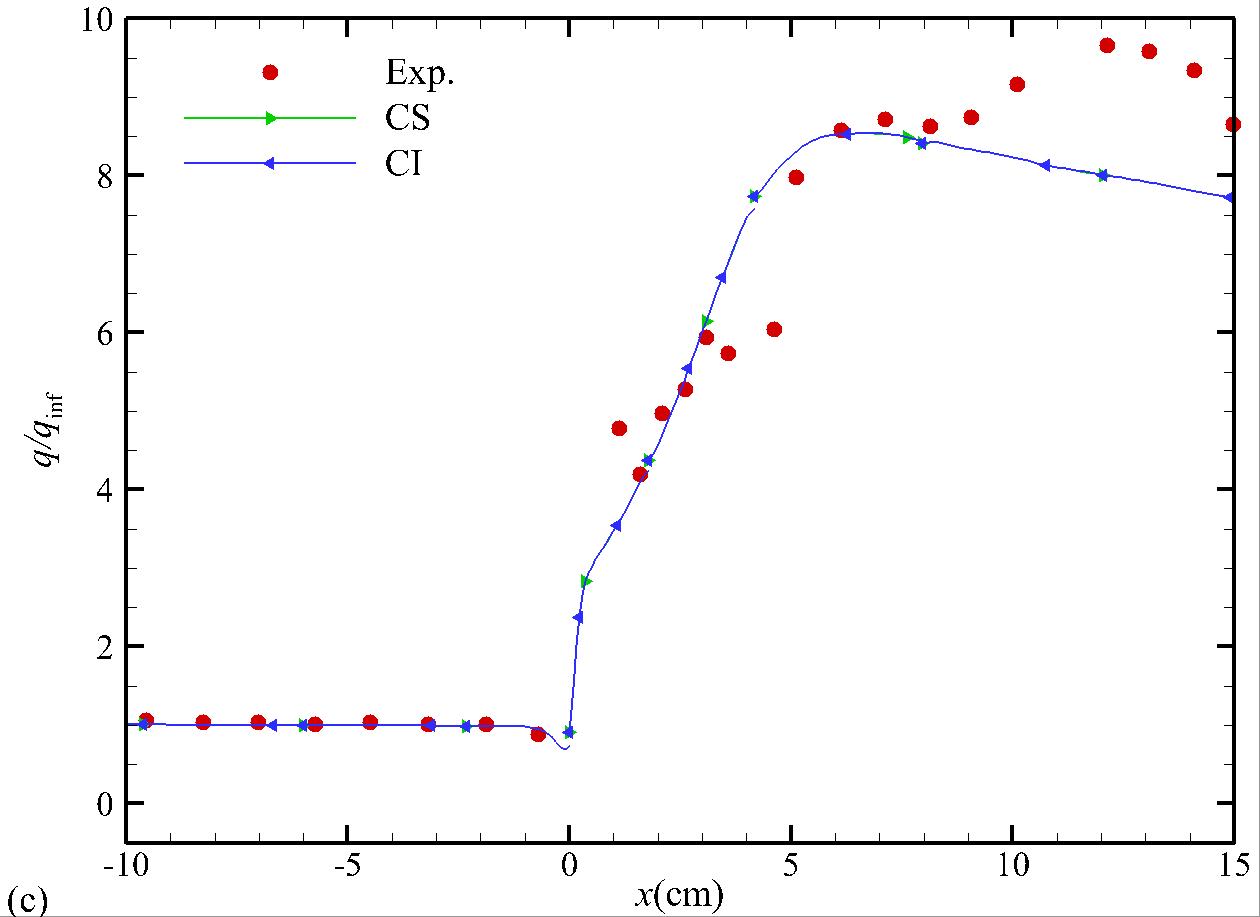}
		% \label{FIG.GSC_grid_a}
	}
	\caption{Comparison of the velocity profile at $x = $-6 cm and normalized surface quantities between the CS, CI, and experimental data for the ASWBLI case:(a) normalized velocity; (b)normalized pressure; (c)normalized wall heat flux.}
	\label{FIG.ASWBLIcontour}
\end{figure}

\begin{figure}[htb!]
	\centering
		\includegraphics[width=0.45\textwidth]{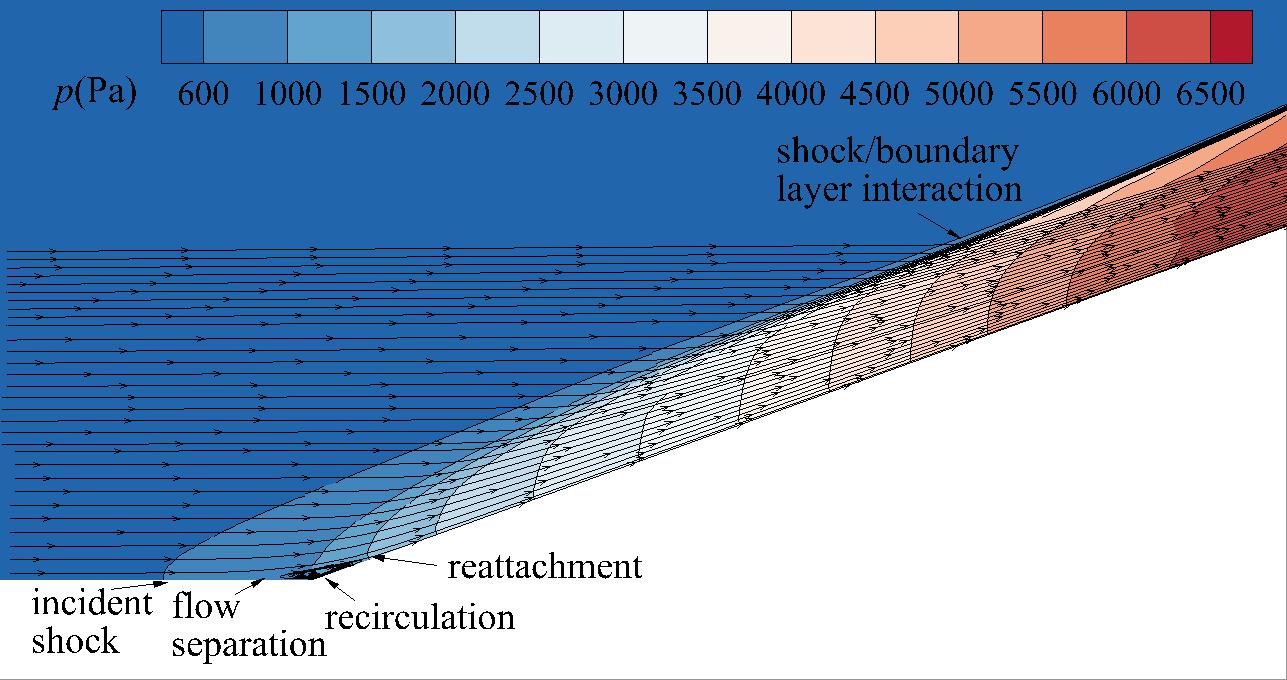}
	\caption{Pressure contour superimposed with streamlines of the ASWBLI case calculated by the CS method.}
	\label{FIG.ASWBLI-pcontour}
\end{figure}

The comparison of results between CS, CI, and experimental data~\cite{Kussoy1989Documentation} are shown in FIG.~\ref{FIG.ASWBLIcontour}. As the residuals have all converged, there are no substantial differences between the computed results of CS and CI among different CFL numbers. In addition, the velocity profile and wall properties from simulation are in good agreement with the experimental results, validating the accuracy of numerical method spatial discretization. The pressure contour superimposed with streamlines computed by the CS method is shown in FIG.~\ref{FIG.ASWBLI-pcontour}. An incident shock wave appears in front of the corner, followed by flow separation and the formation of a recirculation zone, before reattachment occurs. The shock wave and the boundary layer mutually interfere with each other, affecting the positions of both the shock wave and the boundary layer.

\subsection{Hypersonic flows of 3D capsule}

\begin{figure}[htb!]
	\centering
	% \subfigure
	% {
		% \includegraphics[width=.9\textwidth,trim=0 0 0 0,clip]{figure/RAMC_ne.eps}
		\includegraphics[width=.95\textwidth]{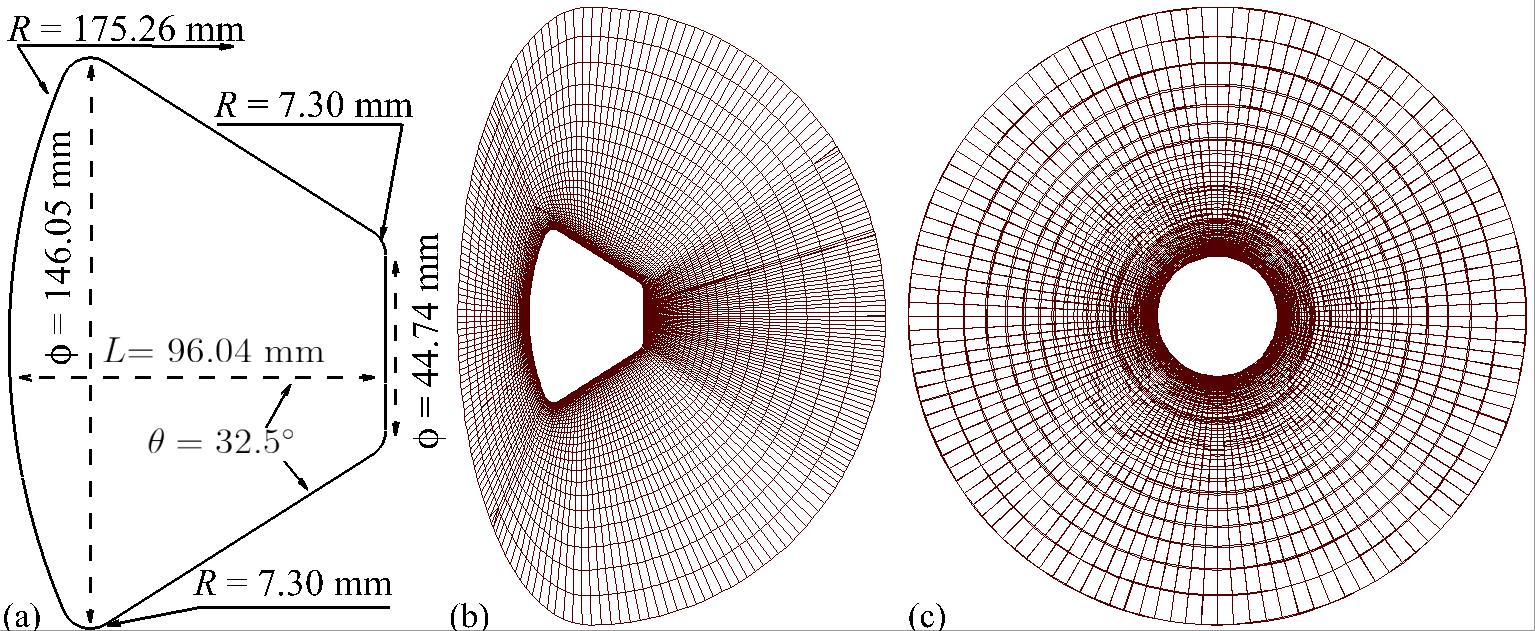}
	% }
	\caption{(a) The geometry of the GSC vehicle; (b) the front view, and (c) the side view of the computational grid of the GSC vehicle.}
	\label{FIG.GSC}
\end{figure}

The well documented generic spherical capsule (GSC)~\cite{maclean2008analysis} is used here to test the component-splitting method on accelerating weakly ionized hypersonic flows. The geometry and computational grid of the GSC model are drawn in FIG.~\ref{FIG.GSC}. The computational body-fitted grid has approximately 0.74 million finite-volume cells with a thickness of the first cell off the wall being $3.29 \times 10^{-6}$ m. Two experimental states are selected for the numerical computation conditions, as shown in Table \ref{Tab.GSCcondition}. The spatial discretization method and chemical kinetic mechanism used to calculate the residual remains consistent with the approach employed in the cylinder case. %The free stream conditions are taken from experimental condition, $u = 2922 \mathrm{m/s}, T=180 \mathrm{K}, \rho = 3.089505 \times 10 ^{-3} \mathrm{kg/m^3}, Y_{N_2} = 0.733, Y_{O_2} = 0.197, Y_{NO} = 0.069$
\begin{table}[htb!]
	\centering
	\caption{Free-stream conditions for numerical simulations of GSC reentry.}
	\label{Tab.GSCcondition}%添加标题 设置标签
	\begin{tabular}{cccccccc}
		\toprule
		%\hline
		case&$u (\mathrm{\frac{m}{s}})$&AOA($^{\circ}$)&$T(\mathrm{K})$&$\rho_{N2}(\mathrm{\frac{kg}{m^3}})$&$\rho_{O2}(\mathrm{\frac{kg}{m^3}})$&$\rho_{NO}(\mathrm{\frac{kg}{m^3}})$&$\rho_{O}(\mathrm{\frac{kg}{m^3}})$\\
		\midrule
		case 1&2922&0.0&180&2.264E{-3}&6.097E{-4}&2.143E{-4}&1.505E{-6}\\
		case 2&2949&28.0&191&2.431E{-3}&6.529E{-4}&2.339E{-4}&1.902E{-6}\\
		%\hline
		\bottomrule
	\end{tabular}
\end{table}

\begin{figure}[htb!]
	\centering
	\subfigure
	{
		\includegraphics[width=0.45\textwidth]{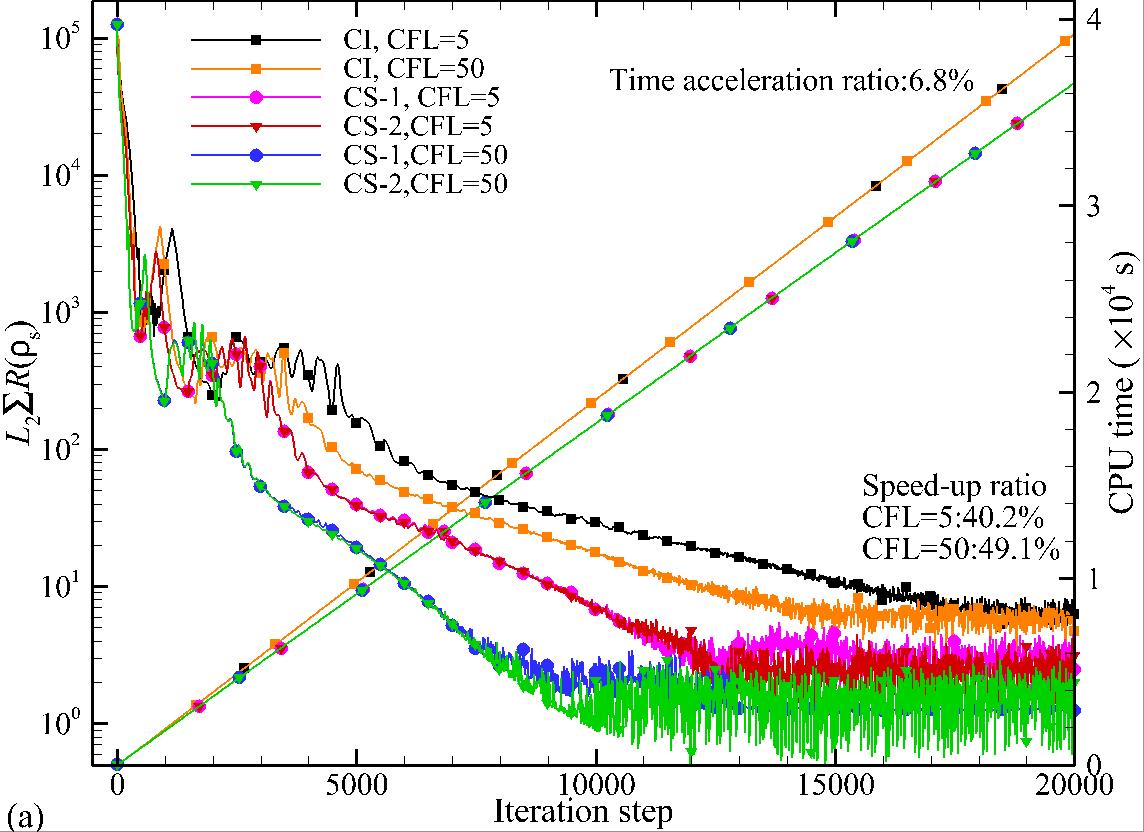}
		% \label{FIG.GSC}
	}
	\subfigure
	{
		\includegraphics[width=0.45\textwidth]{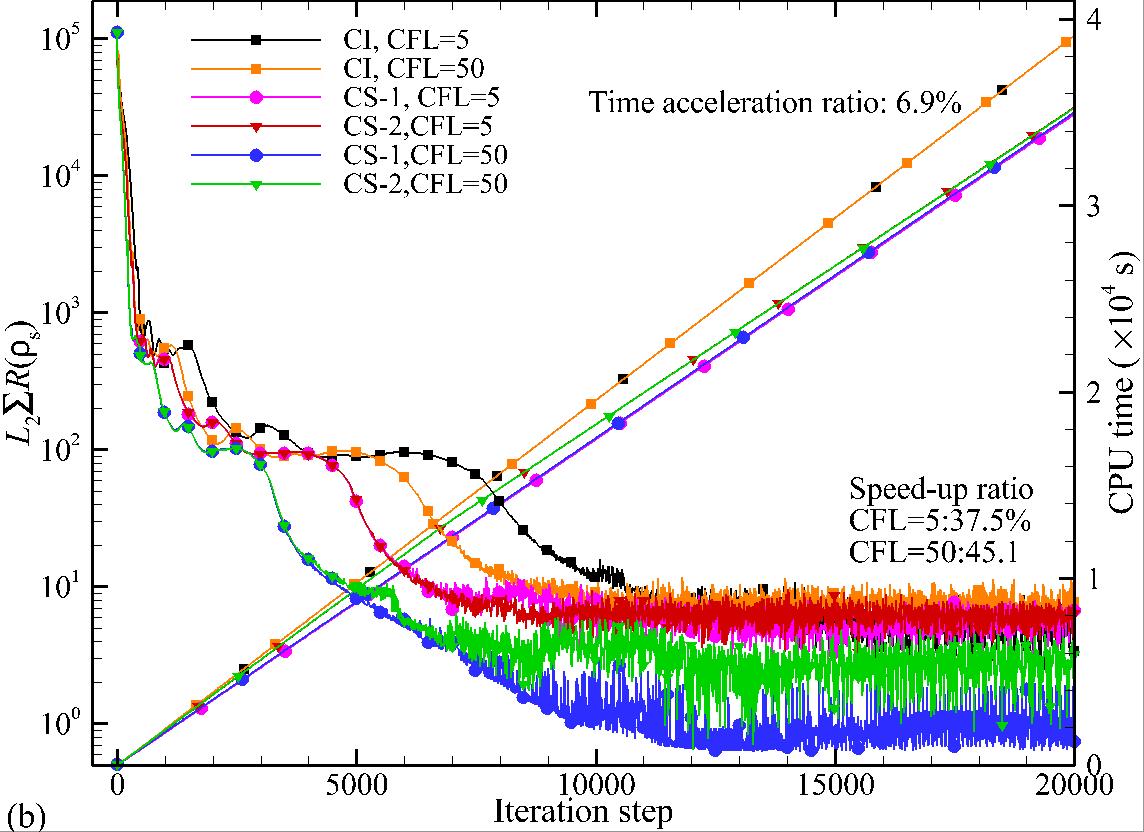}
		% \label{FIG.GSC_grid_a}
	}
	\caption{Comparison of the convergence curve of $L_2$ norm of residual of sum of specie density with respect to iteration step and CPU time between component-splitting scheme and coupled implicit scheme for hypersonic flows of GSC reentry:(a)case 1(b)case 2.}
	\label{FIG.GSCiter}
\end{figure}

The convergence curves of the residual of the species density with respect to iteration step and CPU time for the component-splitting implicit method and the coupled implicit method are shown in FIG.\ref{FIG.GSCiter}. Here, the CS method achieves a reduction in computation time of single-step by approximately 6.9\%. The main acceleration effect is contributed to less iteration steps to achieve convergence, especially in the large CFL number. Consequently, the required number of iterative steps for convergence of case 1 is reduced by 40.2\% and 49.1\% for CFL=5 and CFL=50, respectively. The acceleration ratio of case 2 was similar, and the acceleration effect was more pronounced when a large CFL number was used. In addition, the CS method converges to a lower magnitude of the residual, resulting in higher accuracy.%The maximum stable CFL number of CI for case 1 and case 2 is 1. Under the same small CFL condition, the CS method exhibits a slight acceleration effect compared with the CI method.

\begin{figure}[htb!]
	\centering
	\subfigure
	{
		\includegraphics[width=0.45\textwidth]{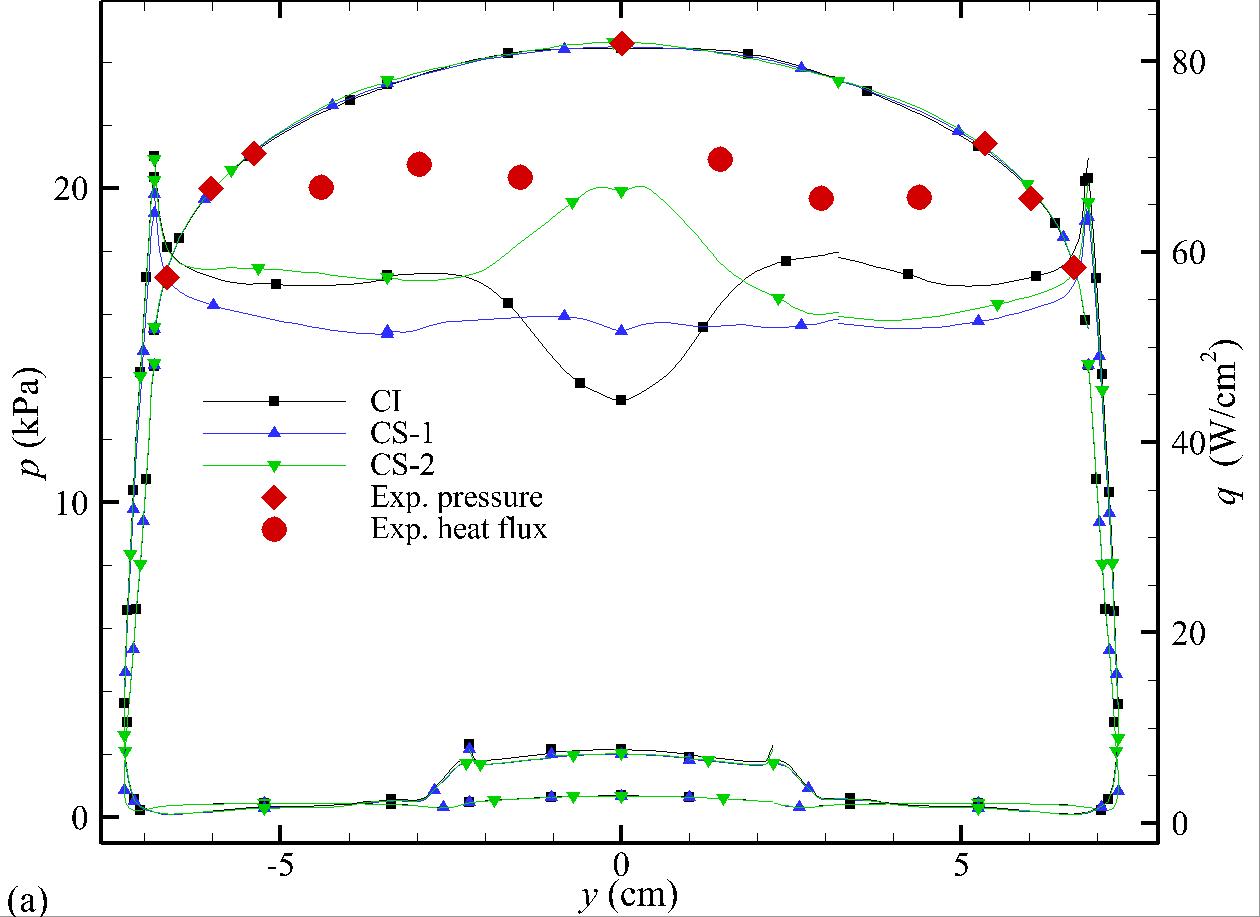}
		% \label{FIG.GSC}
	}
	\subfigure
	{
		\includegraphics[width=0.45\textwidth]{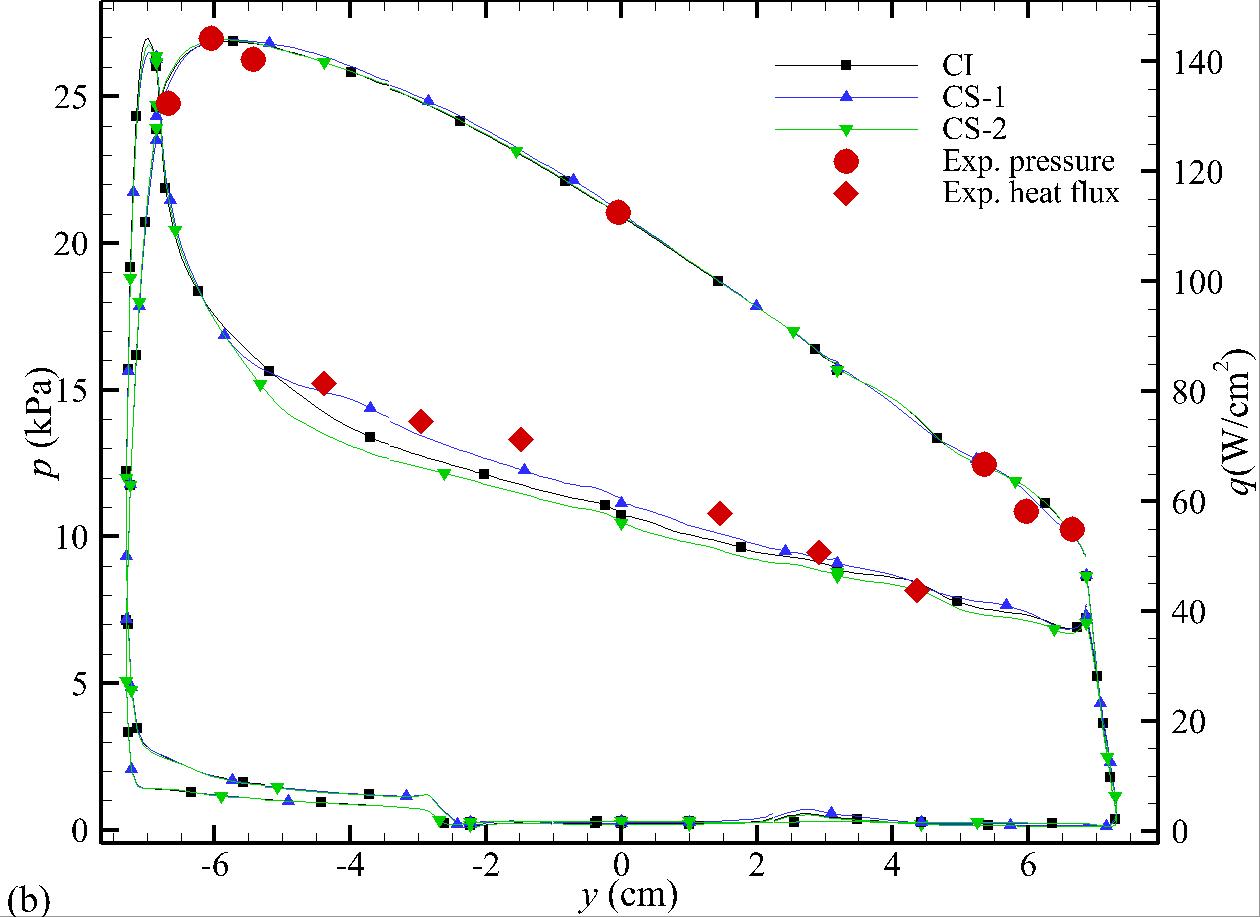}
		% \label{FIG.GSC_grid_a}
	}
	\caption{Comparison of the pressure and heat flux on the surface of GSC reentry between component-splitting scheme and coupled implicit scheme:(a)case 1; (b)case 2.}
	\label{FIG.GSCcompare}
\end{figure}

The required number of iteration steps for the heat flux convergence in the component-splitting method and coupled implicit method are summarized in Table~\ref{Tab.GSCheatflux}. With a convergence criterion of relative error less than 1\%, the component-splitting method accelerates the convergence to heat flux at stagnation points by 50\% and 66.6\% for case 1 and case 2, respectively. The comparison of wall pressure and heat flux between component-splitting method and coupled implicit method with experiment results are shown in FIG.~\ref{FIG.GSCcompare}. The wall pressure computed by the CS and CI methods are in good agreement with experimental results, exhibiting only minor discrepancies. However, differences in wall heat flux are observed between the CS and CI methods. As the residual magnitude of the CS result is smaller, its calculated wall heat flux is more consistent with experimental results~\cite{maclean2008analysis}. The component-splitting method can improve the accuracy of numerical simulations by converging residuals to a lower order of magnitude. 
%The CS method achieves better convergence characeristics of lower residual magnitude while accelerating convergence.

\begin{table}[htb!]
	\centering
	\caption{Comparison of iteration steps to convergence of the wall heat flux at stagnation points in the GSC cases, with a convergence criterion of relative error less than 1\%.}
	\label{Tab.GSCheatflux}%添加标题 设置标签
	\begin{tabular}{ccccc}
		\toprule
		%\hline
		\multirow{2}*{method}&\multicolumn{2}{c}{case 1}&\multicolumn{2}{c}{case 2}\\
		~&CFL=5&CFL=50&CFL=5&CFL=50\\
		\midrule
		CI&8000&6000&8000&6000\\
		CS&6000&3000&5000&2000\\
		%\hline
		\bottomrule
	\end{tabular}
\end{table}

\subsection{Hypersonic flows of winged missile}
\begin{figure}[htb!]
	\centering
	\subfigure
	{
		\includegraphics[height=0.205\textwidth]{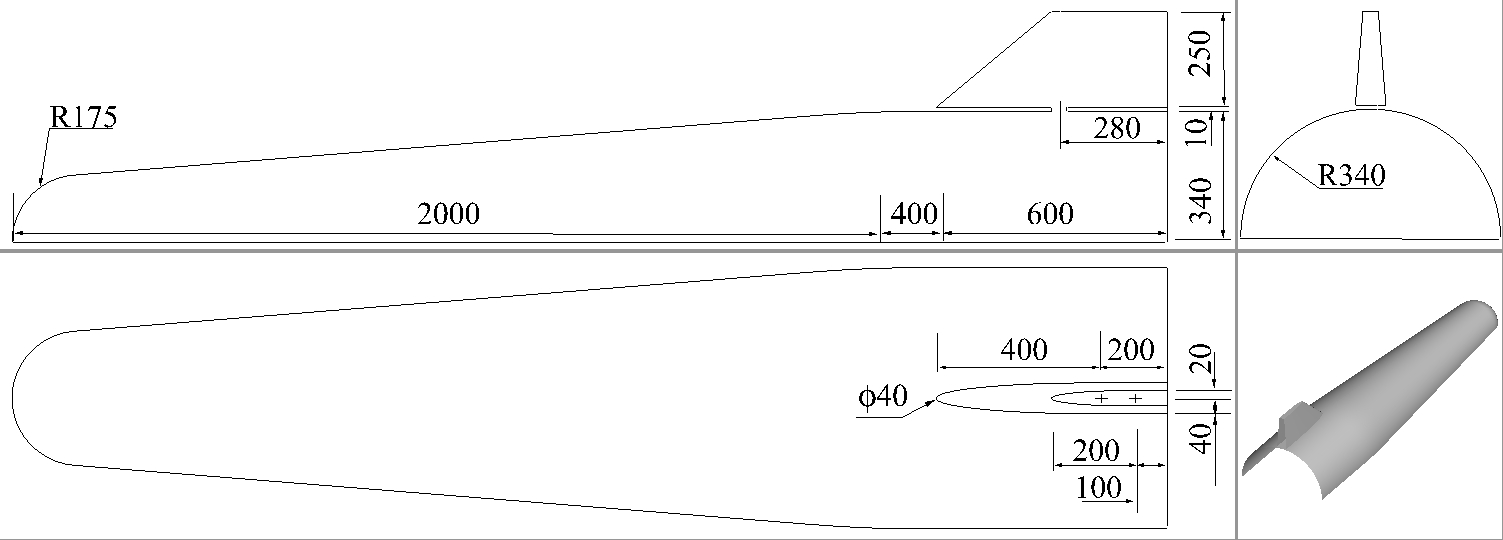}
		% \label{FIG.GSC}
	}
	\subfigure
	{
		\includegraphics[height=0.205\textwidth]{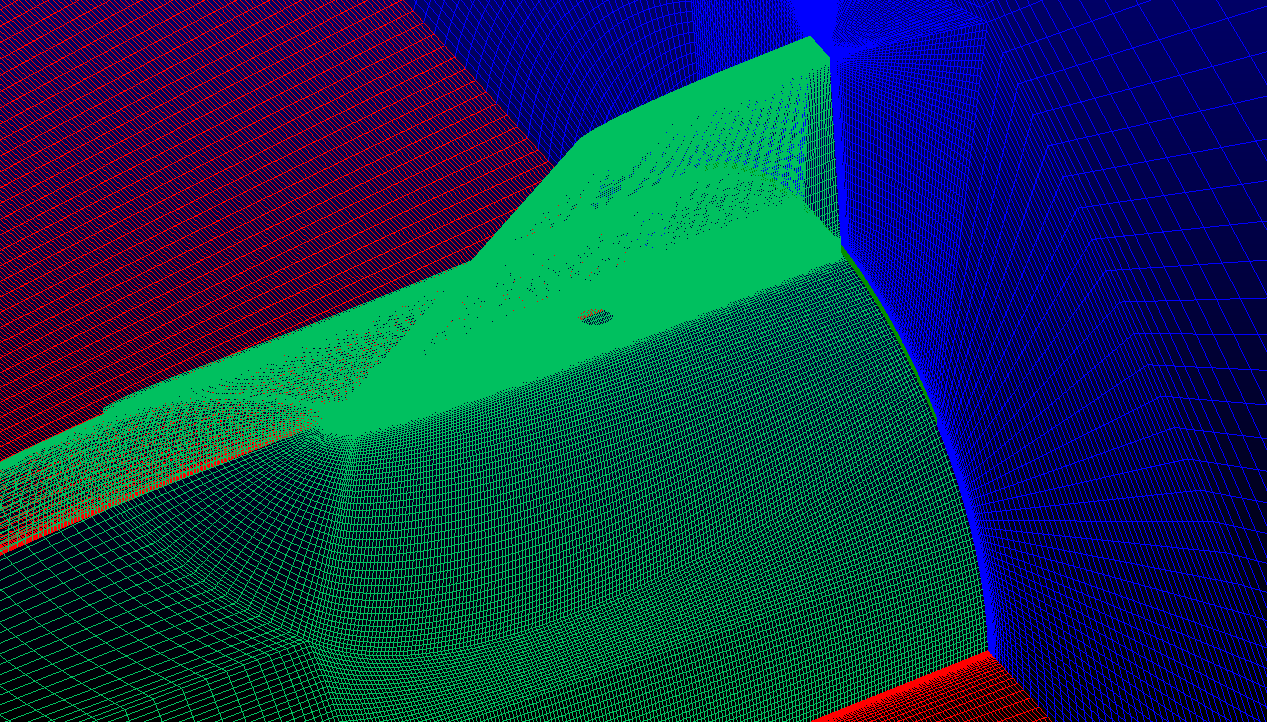}
		% \label{FIG.GSC_grid_a}
	}
	\caption{Geometry and computational grid of the missile (unit:mm).}
	\label{FIG.vehiclegeo}
\end{figure}

To investigate the performance of composent-splitting method in complex flows with large scale grid, numerical test is conducted on a winged missile. The geometry includes the ELECTRE blunt cone~\cite{HAO20161}, a section of cone transition, a section pf cylinder, a tail rudder, and a rotating shaft, as shown in FIG.~\ref{FIG.vehiclegeo}. The computational domain consists of structed body-fitted grid, comprising 4.6 million grid cells in half model. The free-stream condition for numerical simulation is the atmospheric condition at altitude of 60km with a Mach number of 8. The spatial discretization method and chemical kinetic mechanism used to calculate the residual remains consistent with the approach employed in the cylinder case. 

\begin{figure}[htb!]
	\centering
	\subfigure
	{
		\includegraphics[width=0.45\textwidth]{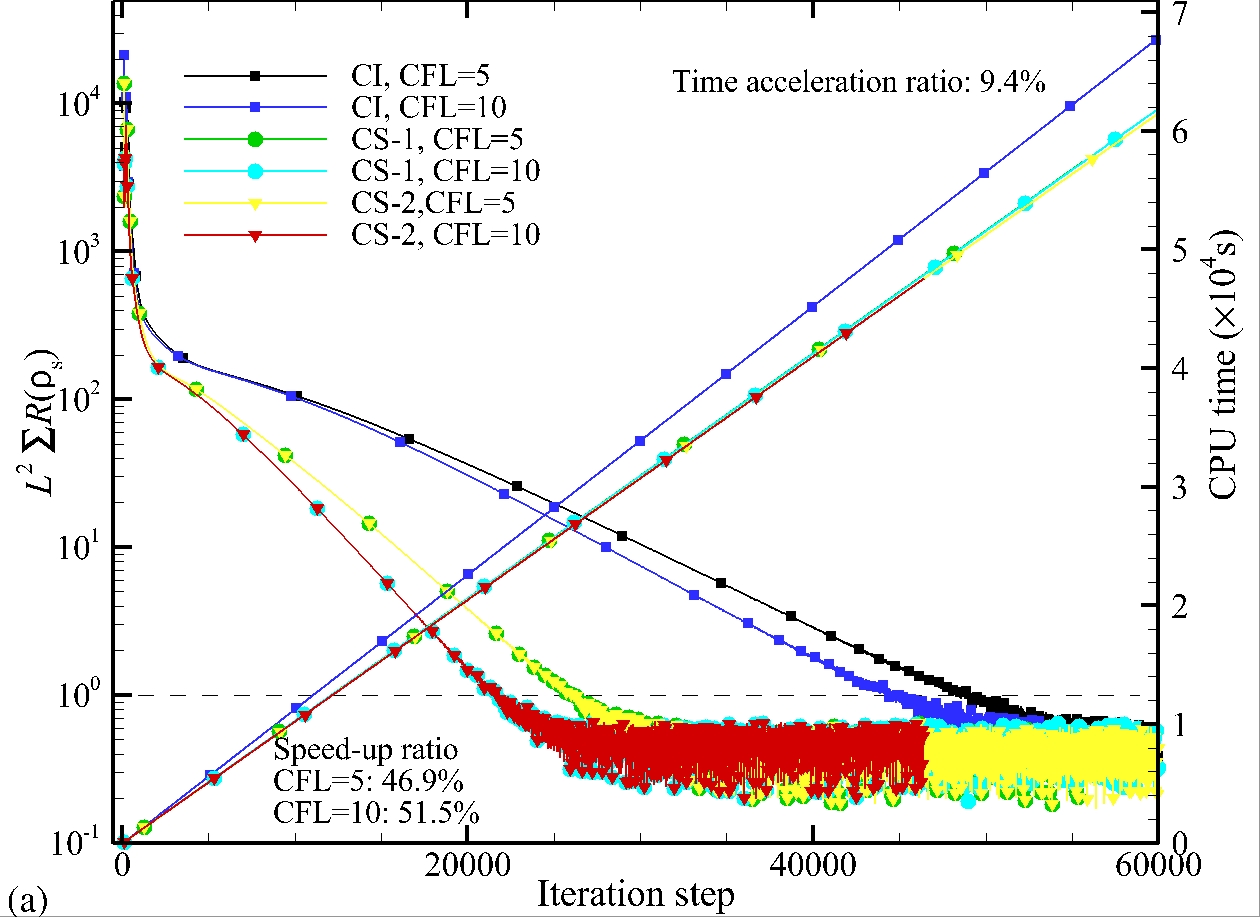}
		% \label{FIG.GSC}
	}
	\subfigure
	{
		\includegraphics[width=0.45\textwidth]{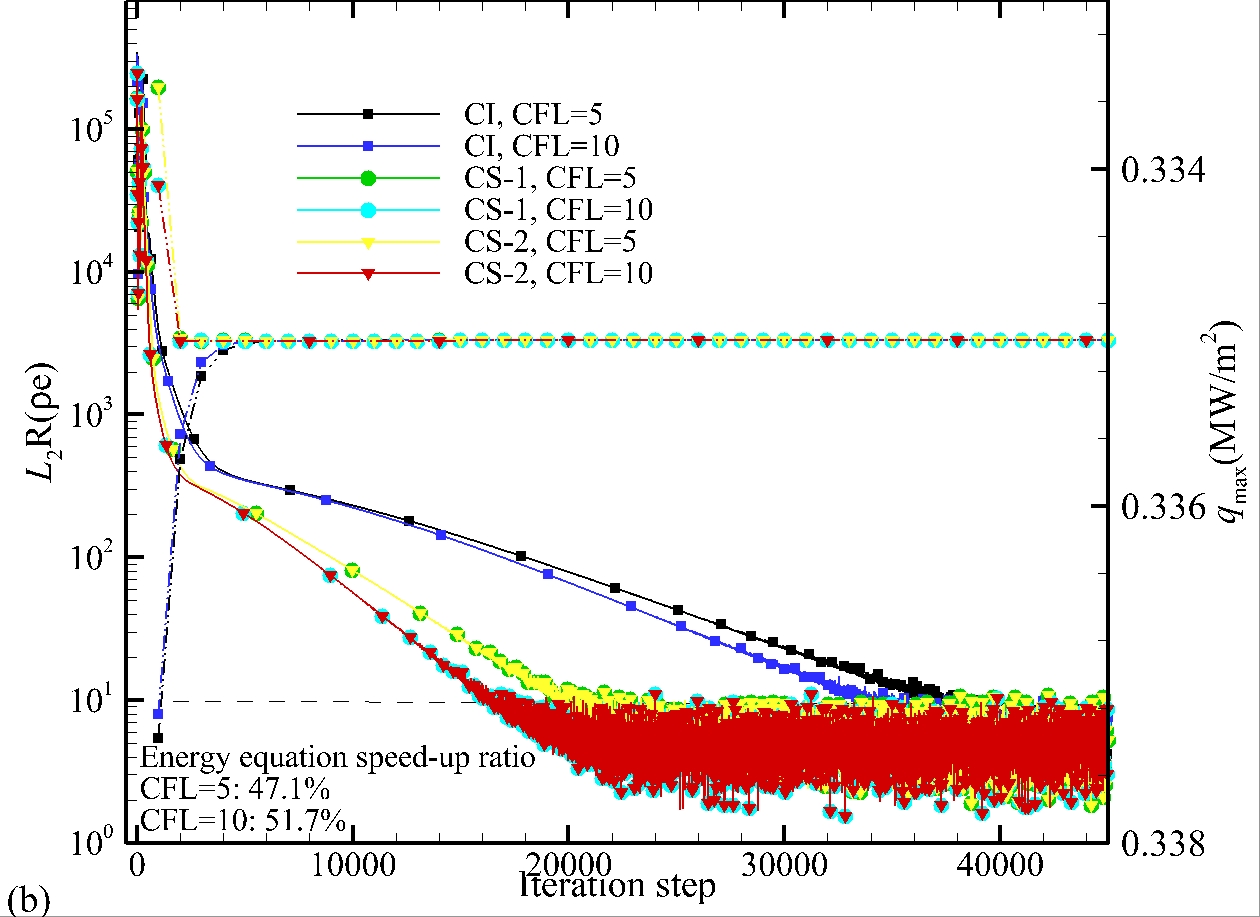}
		% \label{FIG.GSC_grid_a}
	}
	\caption{Comparison of the convergence curve between CS and CI with respect to iteration step in the missile hypersonic case:(a) residual of species density and CPU time; (b)residual of energy and heat flux at stagnation points.}
	\label{FIG.missile-iter}
\end{figure}

\begin{figure}[htb!]
	\centering
	\includegraphics[height=0.205\textwidth]{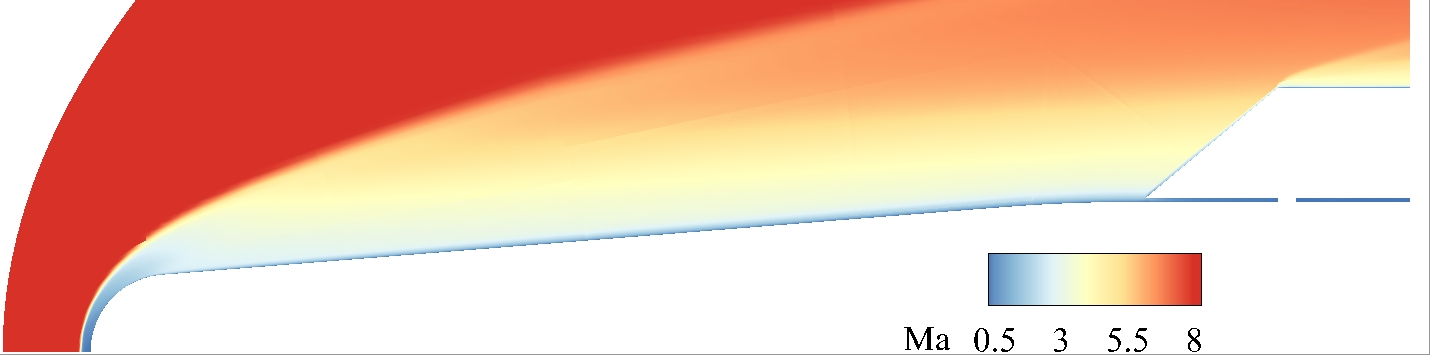}
		% \label{FIG.GSC}
	\caption{Mach number contour of the missile computed by the CS method.}
	\label{FIG.vehiclecontour}
\end{figure}

The convergence curve of heat flux at stagnation points and residual with respet to iteration steps are shown in Fig\ref{FIG.missile-iter}. The component-splitting implicit method achieves substantial convergence acceleation of residual and heat flux. Compared with CI method, the CS method acelerates the convergence of the residual of species density by 46.9\% and 51.5\% for CFL=5 and CFL=10, respectively. The convergence of heat flux and the energy equation also achieves a similar acceleration ratio. In addition, the acceleration ratio of approximatly 9.4\% for a single iteration step of compuation also enhances the computational efficiency.  As shown in FIG.\ref{FIG.vehiclecontour}, the Mach number contour plot was computed using the CS method. The conponent-splitting method can effectively reduce the computational time cost for hypersonic thermochemical nonequilibrium flows.

\section{Conclusion}
In this paper, we propose a component-splitting method to improve the convergence characteristics of the implicit time integration for the multicomponent Navier-Stokes equations. The implicit operator is seprated into the flow parts and component parts due to their distinct eigen system of characteristic decomposition. Using corresponding spectral radius for these type equations in flux splitting to construct implicit operator achieved accelerated convergence. Two consistence corrections have been introduced to avoid the numerical inconsistencies of mass fraction. The first correction normalizes the iterative increment of the conservative variables while ensuring the conservation. The second correction normalizes the mass fraction which compromises the conservation property but achieves better robustness, resulting in higher stable CFL number. Despite compromising the conservation property, the impact of component-splitting method on accuracy is minal when the residual reaches convergence. In addition, the flux jacobian matrix of component equations can be represented with scalar values to achieve higher computational efficiency, avoiding the compuational consumption quadratically increases with the number of species. The time cost of a single implicit iteration is reduced from 33\% to 99 \% with the increasing number of species from 16 to 1024. 

Numerical simulations on hypersonic thermo-chemical nonequalibrium flows are conducted to assess the performance of the component-splitting implicit method in solving multicomponent reactive flows. The convergence acceleration effect has two main aspects: first, the reduction of computation time for each step, and second, the reduction of iteration step to achieve convergence of residual. The component-splitting implicit method also exhibits an accelerated convergence effect in heat flux calculations. The acceleration effects are enhanced with the increasing of CFL number and species number. Lower magnitude of residual is achieved by component-splitting method indicating better accuracy and convergence characteristics.  As a result, the component-splitting method can achieve better convergence characteristics for implicit time integration of multicomponent reactive flows. %A computational time speedup ratio is obtained by 41.9\% for shock/mixinglayer interaction. In all cases, the component-splitting method achieves a smaller residual magnitude and leads to better accuracy.The convergence acceleration effect becomes more pronounced as the CFL number and number of species increased.
%\section*{Appendix}
%An Appendix, if needed, should appear before the acknowledgments.
% \section*{Acknowledgments}
% This work was supported by 111 project on ``Aircraft Complex Flows and the Control'' [grant number B17037].

%% The Appendices part is started with the command \appendix;
%% appendix sections are then done as normal sections
% \section*{Appendix A. Supplementation of the flux jacobian}

%% \section{}
%% \label{}

%% If you have bibdatabase file and want bibtex to generate the
%% bibitems, please use
%%
% \bibliographystyle{elsarticle-num}
%\biboptions{longnamesfirst,angle,semicolon,square,numbers,sort&compress}

%% else use the following coding to input the bibitems directly in the
%% TeX file.

%\begin{thebibliography}{00}

	%% \bibitem[Author(year)]{label}
	%% Text of bibliographic item

%\end{thebibliography}

\begin{acknowledgments}
	This work was supported by 111 project on ``Aircraft Complex Flows and the Control'' [grant number B17037].
	% We wish to acknowledge the support of the author community in using
	% REV\TeX{}, offering suggestions and encouragement, testing new versions,
	% \dots.
\end{acknowledgments}
	
	\section*{Data Availability Statement}
	
	AIP Publishing believes that all datasets underlying the conclusions of the paper should be available to readers. Authors are encouraged to deposit their datasets in publicly available repositories or present them in the main manuscript. All research articles must include a data availability statement stating where the data can be found. In this section, authors should add the respective statement from the chart below based on the availability of data in their paper.

	\bibliography{IterateAcceleration}% Produces the bibliography via BibTeX.
	
	\end{document}